\documentclass[12pt]{article}
\usepackage{amsmath,amsfonts,amssymb}
\usepackage{mathrsfs}
\usepackage{dsfont}
\usepackage[utf8x]{inputenc}
\usepackage[T1]{fontenc}
\usepackage{lmodern} \usepackage{ucs}
\usepackage{fullpage}

\usepackage{graphicx}
\usepackage{caption}
\usepackage{subcaption}

\usepackage{array}
\usepackage{multicol}
\usepackage{multirow}
\usepackage{color}
\usepackage{authblk}

\usepackage{tikz}

\newtheorem{lemma}{Lemma}[section]
\newtheorem{theo}[lemma]{Theorem}
\newtheorem{rmk}[lemma]{Remark}

\newtheorem{defin}[lemma]{Definition}
\newtheorem{coro}[lemma]{Corollary}

\newenvironment{Proof}{\noindent
    \abovedisplayskip = 0.5\abovedisplayskip
    \belowdisplayskip=\abovedisplayskip{\bfseries Proof. }}{\QED\medskip}

    \newenvironment{ProofOf}[1]{\noindent
    \abovedisplayskip = 0.5\abovedisplayskip
    \belowdisplayskip=\abovedisplayskip{\bfseries Proof of  #1. }}{\QED\medskip}

\newcommand{\QED}{\mbox{}\hfill \raisebox{-0.2pt}{\rule{5.6pt}{6pt}\rule{0pt}{0pt}} \medskip\par}

\newcommand{\R}{\mathbb{R}}

\newcommand{\ds}{\displaystyle}
\newcommand{\ud}{\, {\mathrm{d}}}

\title{Particles interacting with a vibrating medium: 
\\
existence of solutions and convergence to the Vlasov--Poisson system}
\author[1]{Stephan De Bi\`evre\thanks{{\tt stephan.de-bievre@univ-lille1.fr}}}
\author[2]{Thierry~Goudon\thanks{ {\tt thierry.goudon@inria.fr}}}
\author[2]{Arthur Vavasseur\thanks{ {\tt arthur.vavasseur@unice.fr}}}

\affil[1]{\small Univ. Lille, CNRS, UMR 8524 - Laboratoire Paul Painlev\'e, F-59000 Lille, France\thanks{S.D.B. acknowledges the support of the Labex CEMPI (ANR-11-LABX-0007-01).
}.

\&
Equipe-Projet MEPHYSTO,
Centre de Recherche INRIA Futurs,

Parc Scientifique de la Haute Borne, 40, avenue Halley B.P. 70478,\\
F-59658 Villeneuve d'Ascq cedex, France.}
\affil[2]{\small Inria,  Sophia Antipolis M\'editerran\'ee Research Centre, Project COFFEE

\& Univ. Nice Sophia Antipolis, CNRS, Labo. J. A. Dieudonn\'e, UMR 7351 

Parc Valrose, F-06108 Nice, France}

\begin{document}
\maketitle

\abstract{We are interested in a kinetic equation intended to describe the interaction of particles with their environment.
The environment is modeled by a collection of local vibrational degrees of freedom.
We establish the existence of weak solutions for a wide class of initial data  and  external forces. 
We also identify a relevant regime which allows us to derive, quite surprisingly, the attractive Vlasov--Poisson system from the coupled Vlasov-Wave equations.}

\vspace*{.5cm}
{\small
Inelastic Lorentz gas. \\

\noindent{\bf Math.~Subject Classification.} 
82C70, 
70F45, 
37K05, 
74A25. 
}

\section{Introduction}

In \cite{BdB}, L. Bruneau and S. De Bi\`evre introduced a mathematical model intended to describe the 
interaction of a classical particle with its environment.
The environment is modeled by a vibrating scalar field, and the dynamics is governed by energy exchanges between the particle and the field, 
embodied into a Hamiltonian structure.
To be more specific on the model in \cite{BdB}, let us denote by $q(t)\in\mathbb R^d$ the position occupied by the particle at time $t$.
The environment
is represented by a field $(t,x,y)\in \mathbb R\times \mathbb R^d\times \mathbb R^n\mapsto \Psi(t,x,y)\in \mathbb R$:
it can be thought of as an infinite set of $n$-dimensional membranes, one for each $x\in\R^d$. The displacement of the membrane positioned at $x\in\R^d$ is given by $y\in\R^n\mapsto \psi(t,x, y)\in\R$.  
 The coupling  is realized by means of  form factor functions $x\mapsto \sigma_1(x)$ and $y\mapsto \sigma_2(y)$, which are supposed to be 
 non-negative,
 infinitely smooth, radially symmetric  and compactly supported. 
 Therefore, the dynamic is described by the following set of differential equations
\begin{equation}\label{sdb} \left\lbrace \begin{array}{l} 
\ddot{q}(t) = -\nabla V(q(t)) - \ds\int_{\mathbb R^d\times \mathbb R^n} \sigma_1 (q(t)-z)\ \sigma_2(y) \ \nabla_x \Psi(t,z,y) \ud y \ud z,
 \\
 [.4cm]
\partial_{tt}^2 \Psi (t,x,y) - c^2 \Delta_y \Psi(t,x,y) = - \sigma_2(y)\sigma_1 (x-q(t)) , \qquad x \in \mathbb R^d ,\  y \in \mathbb R^n. 
\end{array} \right. \end{equation}
In \eqref{sdb}, $c>0$ stands for the wave speed in the transverse direction, while $q\in\mathbb R^d\mapsto V(q)\in\mathbb R$
is a time-independent external potential the particle is subjected to.
In \cite{BdB}, the well-posedness
theory for \eqref{sdb} is investigated, but the main issue addressed there is the large time behavior of the system. 
It  is shown that the system exhibits  dissipative features: under certain circumstances
(roughly speaking, $n=3$ and $c$ large enough) and for a large class of finite energy initial conditions
the particle energy is evacuated in the membranes, and the environment acts with  a friction force on the particle.
Accordingly, the asymptotic behavior of the particle for large times can be characterized
depending on the external force: if $V=0$, the particle stops exponentially fast, when $V$ is a confining potential
 with a minimiser $q_0$, then the particle stops at the location $q_0$, and for $V(q)=-F\cdot q$, a limiting velocity $V_F$ can be identified.
 \\
 
 Since then, a series of works has been devoted to further investigation of the asymptotic properties of a family of related models.
 We refer the reader to  \cite{AdBLP,dBLP,dBP,dBPS,LPdB,dBS} for thorough numerical experiments
 and analytical studies, that use random walks arguments in particular.
 The model can be seen as a variation on the  Lorentz gas model where one is interested in  
 the free motion of a
single point particle in a system of obstacles distributed on a certain lattice.
 We refer the reader to \cite{BBS, CG,Gal, FG_Lo, MS} for results and 
 recent overviews on the Lorentz gas problem.
 Instead 
of dealing with periodically or randomly distributed hard scatterers as in the Lorentz gas model, here
the particle interacts with a vibrational environment, that create the ``soft'' potential $\Phi$.
 The asymptotic analysis of the behavior of a particle subjected to an oscillating 
 potential is a further related problem that is also worth  mentioning \cite{FH,GR, KP,  PV}.
 \\
 
We wish to revisit the model of \cite{BdB} in the framework of kinetic equations.
Instead of considering a single particle described by its position $t\mapsto q(t)$, we work with the particle distribution function in phase space 
$f(t,x,v)\geq 0$, with $x\in\mathbb R^d$, $v\in\mathbb R^d$, the position and velocity variables respectively.
This quantity obeys the following Vlasov equation
\begin{equation}\label{kin}
\partial_t f + v\cdot\nabla_x f - \nabla_x (V+\Phi)\cdot \nabla_v f =0,  
\qquad t\geq 0,\ x\in \mathbb R^d,\ v\in \mathbb R^d.
\end{equation}
In \eqref{kin}, $V$ stands for the external potential, while $\Phi$ is the self-consistent potential
describing the interaction with the environment.
It is defined by the convolution formula
\begin{equation}\label{pot}
\Phi(t,x)=\ds \int_{\mathbb R^d \times \mathbb R^n} \Psi(t,z,y) \sigma_1(x-z) \sigma_2(y) \ud y\ud z,  \qquad t\geq 0,\ x\in \mathbb R^d
\end{equation}
where the vibrating field 
$\Psi$ is driven by the following wave equation
\begin{equation} 
\label{modele_milieu}
\left\lbrace \begin{array}{l} 
\big(\partial_{tt}^2 \Psi -c^2 \Delta_y \Psi\big)(t,x,y) = -
\sigma_2(y)\ \ds\int_{\mathbb R^d} \sigma_1(x-z) \rho(t,z)\ud z,\
 t\geq 0,\ x\in \mathbb R^d,\ y\in \mathbb R^n,
 \\
 [.4cm]
 \rho(t,x)=\ds\int_{\mathbb R^d} f(t,x,v)\ud v.\end{array}\right.\end{equation}
The system is completed by initial data
\begin{equation} 
\label{CI}
f(0,x,v) = f_0(x,v),
\qquad
\Psi(0,x,y) = \Psi_0 (x,y) ,\qquad 
\partial_t \Psi(0,x,y) = \Psi_1(x,y). 
\end{equation}
A possible interpretation of the kinetic equation \eqref{kin}
consists in considering the model \eqref{sdb} for a set of $N\gg 1$ particles.
The definition of the self--consistent potential has to be adapted since all the particles interact with the environment, namely we have, for $j\in\{1,...,N\}$
 \[ \left\lbrace \begin{array}{l} 
\ddot{q}_j(t) = -\nabla V(q_j(t)) - 
\ds\int_{\mathbb R^d\times \mathbb R^n} \sigma_1 (q_j(t)-z)\ \sigma_2(y) \ \nabla_x \Psi(t,z,y) \ud y \ud z,
 \\
 [.4cm]
\partial_{tt}^2 \Psi (t,x,y) - c^2 \Delta_y \Psi(t,x,y) = - \sigma_2(y)\ \ds\sum_{k=1}^N\sigma_1 (x-q_k(t)). 
\end{array} \right. \]
Note that such a many-particle system is not considered in \cite{BdB}. It is very likely that its asymptotic behavior is much more complicated than with a single particle because, even if the particles do not interact directly, they do so indirectly via their interaction with the membranes. 
If we now adopt the mean--field rescaling in which $\Phi\rightarrow \frac1N\Phi$, then 
\eqref{kin} can be obtained as the limit as $N$ goes to $\infty$
for the empirical measure $f_N(t,x,v)=\frac1N \sum_{k=1}^N \delta(x=q_k(t),v=\dot{q}_k(t))$ of the $N-$particle system,
 assuming the convergence of the initial state $f_N(0,x,v)\rightarrow f_0(x,v)$ in some suitable sense.
 Such a statement can be rephrased 
in terms of the convergence of the joint distribution of the 
$N$--particle system. This issue will be discussed elsewhere \cite{AVPhD} and 
we refer the reader to the lecture notes \cite{FG} and to \cite{GMR}
for further information on the mean--field regimes in statistical physics.
\\

In this paper we wish to analyse several aspects of the Vlasov-Wave system \eqref{kin}--\eqref{CI}.
We warn the reader that, despite the similarities in terminology, the model considered here is very different, both mathematically and physically, from the one dealt with in \cite{BGP}, which is a simplified version of the 
Vlasov--Maxwell system.  It is indeed crucial to understand that the wave equation in this paper is set with variables \emph{transverse} to the physical space: the waves do not propagate at all in the space where the particles move.  This leads to very different physical effects; we refer to~\cite{BdB} and references therein for more details on this matter. 
We add that this paper is less ambitious than \cite{BdB}, since we do not discuss here the large time behavior of the solutions, only their global 
existence.
As mentioned above, since we are dealing with many particles, it is very likely that the question cannot be handled in the same terms as in \cite{BdB}, and 
that the kinetic model inherits the same technical and conceptual difficulties already mentioned for $N>1$ particles.
We only mention that a particular stationary solution (with $f$ integrable) has been exhibited in  \cite{AGV}, and that this solution is shown  to be linearly stable.
\\

The paper is organized as follows.
Section~\ref{s:prelim} contains a preliminary and largely informal discussion to set up notation and to establish some estimates on the interaction potential needed in the bulk of the paper. 
Section~\ref{Cau} establishes the well--posedness of the problem  \eqref{kin}--\eqref{CI} (Theorem~\ref{th1}).
We consider a large class of initial data and external potentials with functional arguments which are reminiscient of Dobrushin's analysis
of the Vlasov equation \cite{Dob}. 
Section~\ref{As} is devoted to
asymptotic issues which allow us to connect 
\eqref{kin}--\eqref{CI} to Vlasov equations with an \emph{attractive} self--consistent 
potential.
In particular, up to a suitable rescaling of the form function $\sigma_1$, we can derive this way the attractive Vlasov--Poisson system. This is quite surprising 
and unexpected in view of the very different 
physical motivation of the models.


\section{Preliminary discussion}\label{s:prelim}
Throughout the paper, we make the following assumptions on the model parameters and on the initial conditions. First, on the coupling functions $\sigma_1, \sigma_2$, we impose:
\begin{equation}\tag{{\bf{H1}}}\label{H1}
\left\{
\begin{array}{l}
\sigma_1\in C^\infty_c(\mathbb R^d,\mathbb R), \ \sigma_2\in C^\infty_c(\mathbb R^n,\mathbb R),
\\
\sigma_1(x)\geq 0,\ \sigma_2(y)\geq 0\  \text{for any $x\in\mathbb R^d$, $y\in\mathbb R^n$},
\\
\text{$\sigma_1,\sigma_2$ are radially symmetric}.
\end{array}
\right.\end{equation}
We require that the external potential fulfills 
\begin{equation}\tag{{\bf{H2}}}\label{H2}
\left\{\begin{array}{l}V\in W^{2,\infty}_{\mathrm{loc}}(\mathbb R^d),
\\
\text{and there exists $C\geq0$ such that   $V(x) \geq -C(1+|x|^2)$ for any $x\in\mathbb R^d$.} 
\end{array}\right.\end{equation}
This is a rather standard and natural assumption. Note that it ensures global existence when $\sigma_1=0=\sigma_2$: it then 
implies that the external potential cannot drive the particle to infinity in finite time.
For the initial condition of the vibrating environment, we shall assume 
\begin{equation}\tag{{\bf{H3}}}\label{H3}
\Psi_0, \Psi_1\in 
 L^2(\mathbb R^d\times\mathbb R^n).
\end{equation}
For the 
 initial particle distribution function, we naturally assume
\begin{equation}\tag{{\bf{H4}}}\label{H5}
f_0\geq 0,\qquad f_0\in L^1(\mathbb R^d\times\mathbb R^d).
\end{equation}
For energy considerations, it is also relevant to suppose
\begin{equation}\tag{{\bf{H5}}}\label{H6}
\nabla_y\Psi_0\in 
 L^2(\mathbb R^d\times\mathbb R^n)\quad  \textrm{ and } \quad \Big((x,v)\mapsto (V(x)+|v|^2)f_0(x,v)\Big)\in L^1(\mathbb R^d\times\mathbb R^d).
\end{equation}
This means that the initial state has finite mass, potential and kinetic energy.

Our goal in this section is to rewrite the equations of the coupled system~\eqref{kin}-\eqref{CI} in an equivalent manner, more suitable for our subsequent analysis.  The discussion will be informal, with all computations done for sufficiently smooth solutions. The proper functional framework will be provided in the next section. First, we note that 
it is clear that \eqref{kin} preserves the total mass of the particles
\[
\ds\frac{\ud}{\ud t} \ds\int_{\mathbb R^d\times \mathbb R^d} f(t,x,v)\ud v\ud x=0.
\] 
In fact, since the field $(v,\nabla_x V+\nabla_x\Phi)$ is divergence--free (with respect to the phase variables $(x,v)$), 
 any $L^p$ norm of the density $f$ is conserved, $1\leq p\leq \infty$.
Furthermore, the PDEs system \eqref{kin}--\eqref{modele_milieu} inherits from the Hamiltonian nature of the original equations of motion~\eqref{sdb} the following 
easily checked energy conservation property:
\[\begin{array}{l}
\ds\frac{\ud}{\ud t}\left\{ \frac{1}{2} \int_{\mathbb R^d\times \mathbb R^n} \left| \partial_t \Psi(t,x,y) \right|^2 \ud x  \ud y 
+ \frac{c^2}{2} \int_{\mathbb R^d\times \mathbb R^n} \left| \ds\nabla_y \Psi(t,x,y) \right|^2 \ud x  \ud y 
\right.
\\[.4cm]\left.
\qquad\qquad\qquad\qquad\qquad + \ds\int_{\mathbb R^d\times \mathbb R^d} f(t,x,v) \left( \frac{|v|^2}{2} + V(x) + \Phi(t,x) \right) \ud v  \ud x
\right\}=0. \end{array} \]
As a matter of fact the energy remains finite when the full set of assumptions \eqref{H1}--\eqref{H6} holds.


For the Vlasov--Poisson equation it is well known that the potential can be expressed by means of a convolution formula.
Similarly here, the self-consistent potential $\Phi$ can be computed explicitly as the image  
of a certain linear operator acting on the macroscopic density $\rho(t,x)=\int_{\mathbb R^d}f(t,x,v)\ud v$; this follows from the fact that the linear wave equation~\eqref{modele_milieu} can be solved explicitly as the sum of the solution of the homogeneous wave equation with the correct initial conditions plus the retarded solution of the inhomogeneous wave equation. To see how this works,
we introduce
\[ 
t\mapsto p(t) = \ds\frac{1}{(2\pi)^n} \int_{\mathbb R^n} \frac{\sin(c|\xi|t)}{c |\xi|} \left| \widehat{\sigma_2}(\xi) \right|^2 \ud\xi\]
and
\begin{equation}
\Phi_0 (t,x) = \frac{1}{(2\pi)^n} 
\ds\int_{\mathbb R^n}\ds\int_{\mathbb R^d}
 \sigma_1 (x-z)
\left( \widehat{\Psi_0}(z,\xi) \cos(c|\xi|t) + \widehat{\Psi_1}(z,\xi) \frac{\sin(c|\xi|t)}{c|\xi|} \right) \ \widehat{\sigma_2}(\xi) \ud z\ud \xi \label{eq:Phizero}
\end{equation}
where the symbol $\, \widehat{\cdot}\, $ stands for the Fourier transform with respect to the variable  $y\in\mathbb R^n$. Note that $\Phi_0$ is the solution of the homogeneous wave equation with the given initial conditions for $\Psi$. 
Finally, we define the operator $\mathcal L$ which associates to a distribution function 
$f:(0,\infty)\times \mathbb R^d\times\mathbb R^d\rightarrow \mathbb R$ the quantity
\begin{equation}\label{eq:Loperator} 
\mathcal{L}(f)(t,x) = \int_0^t p(t-s)\left(\int_{\mathbb R^d}\Sigma(x-z) \rho(s,z)\ud z\right)\ud s,
\end{equation}
where
\[\rho(t,x)=\ds\int_{\mathbb R^d}f(t,x,v)\ud v,\qquad \Sigma=\sigma_1\underset{x}{\ast}\sigma_1. \]
We can then check that 
 the pair $(f,\Psi)$ is a solution of  \eqref{kin}--\eqref{modele_milieu} iff $f$ satisfies 
 \begin{equation} 
\label{simplification}
\left\lbrace \begin{array}{l}
\partial_t f + v\cdot\nabla_x f = \nabla_v f \cdot \nabla_x \left(V + \Phi_0 - \mathcal{L}(f) \right) \\
f(0,x,v)= f_0(x,v) \end{array} \right. 
\end{equation}
and  $\Psi$ is the unique solution of \eqref{modele_milieu}.

We sketch the computation, which is instructive. 
Let  $(f,\Psi)$ be a  solution of \eqref{kin}--\eqref{modele_milieu}. Applying the  Fourier transform with respect to the  variable $y$ 
we find \[\left\lbrace \begin{array}{l} 
(\partial_t^2 + c^2 |\xi|^2) \widehat{\Psi}(t,x,\xi) = - (\rho(t,\cdot) \underset{x}\ast \sigma_1 )(x)\ \widehat{\sigma_2}(\xi), \\[.4cm]
\widehat{ \Psi}(0,x,\xi) = \widehat{\Psi_0} (x,\xi) \qquad
\partial_t \widehat{\Psi}(0,x,\xi) = \widehat{\Psi_1}(x,\xi) .
\end{array} \right.\]
The solution reads
\begin{equation}
\label{Fourrier}\begin{array}{lll}
\widehat{\Psi}(t,x,\xi) &=& - \ds\int_0^t (\rho(t-s,\cdot) \ast \sigma_1 )(x)\  \widehat{\sigma_2}(\xi) \frac{\sin(c s|\xi|)}{c|\xi|} \ud s
\\[.4cm]&&\qquad\qquad+\ds \widehat{\Psi_0}(x,\xi) \cos(c|\xi|t) + \widehat{\Psi_1}(x,\xi) \frac{\sin(c|\xi|t)}{c|\xi|}. 
\end{array}\end{equation}
To compute $\Phi$ in~\eqref{pot}, we use Plancherel's equality:
\[ \Phi(t,x) \begin{array}[t]{l} \ds   = \int_{\mathbb R^d \times \mathbb R^n} \Psi(t,z,y) \sigma_1(x-z) \sigma_2(y) \ud y\ud z \\[.4cm] \ds 
=\frac{1}{(2\pi)^n}  \int_{\mathbb R^d \times \mathbb R^n} \widehat{\Psi}(t,z,\xi) \sigma_1(x-z) \widehat{\sigma_2}(\xi) \ud \xi \ud z  \\[.4cm] \ds 
=-\left((\sigma_1 \ast \sigma_1) \ast
\ds\int_0^t \Big(\rho(t-s,\cdot) \int_{\mathbb R^n} \frac{\sin(c s|\xi|)}{c|\xi|} \frac{\left| \widehat{\sigma_2}(\xi) \right|^2}{(2\pi)^n} \ud \xi \Big)\ud s\right)(x)
\\[.4cm]
\qquad+ \ds\frac{1}{(2\pi)^n}
\left(\sigma_1 \ast\ds \int_{\mathbb R^n} \left(\widehat{\Psi_0}(\cdot,\xi) \cos(c|\xi|t) + \widehat{\Psi_1}(\cdot,\xi) \ds\frac{\sin(c|\xi|t)}{c|\xi|} \right)\ \widehat{\sigma_2}(\xi) \ud\xi\right)(x) \\ [.4cm]\ds 
=- \mathcal{L}(f)(t,x) + \Phi_0(t,x).
\end{array}\]
Inserting this relation into   \eqref{kin}, we arrive at \eqref{simplification}. Conversely, let $f$ be a  solution of \eqref{simplification} and let  $\Psi$  be the unique solution of \eqref{modele_milieu}.
The same 
computation then shows that $\Phi$ in~\eqref{pot} is given by
 $\Phi = \Phi_0 - \mathcal{L}(f)$. Therefore $f$ satisfies \eqref{kin}.
\\

The operator ${\mathcal L}$ in~\eqref{eq:Loperator} plays a crucial role in our further analysis. Its precise definition on an appropriate functional space and its basic continuity properties are given in the following Lemma.
\begin{lemma}[Estimates on the interaction  potential]\label{l:estL}
For any  $0<T<\infty$, the following properties hold:
\begin{itemize}
\item[i)]    $\mathcal{L}$ belongs to the space $\mathcal{A}_T$ of  continuous operators on 
$C \big( [0,T];\big( W^{1,\infty}(\mathbb R^d\times \mathbb R^d) \big)' \big)$ with values in 
$C\big([0,T];W^{2,\infty}(\mathbb R^d)\big)$. 
Its norm is evaluated 
as follows:
\[ |\!|\!| \mathcal{L} |\!|\!|_{\mathcal{A}_T} \leq \| \sigma_1 \|^2_{W^{3,2}(\mathbb R^d)} 
\| \sigma_2 \|^2_{L^2(\mathbb R^n)} \ \ds\frac{T^2}{2}; \]
\item[ii)]  $\mathcal{L}$ belongs to the space $\mathcal{B}_T$ of  continuous operators on 
$C \big( [0,T];\big( W^{1,\infty}(\mathbb R^d\times \mathbb R^d) \big)' \big)$ with values in 
$C^1\big([0,T];L^\infty(\mathbb R^d)\big)$. Its norm is evaluated 
as follows:
\[ |\!|\!|  \mathcal{L} |\!|\!| _{\mathcal{B}_T} \leq \| \sigma_1 \|^2_{W^{1,2}(\mathbb R^d)} \|\sigma_2\|^2_{L^2(\mathbb R^n)} \ \Big(T+ \ds\frac{T^2}{2}\Big); \]
\item[iii)] $ \Phi_0 $ satisfies
\[
\| \Phi_0(t,\cdot) \|_{W^{2,\infty}(\mathbb R^d)} \leq 
\| \sigma_1 \|_{W^{2,2}(\mathbb R^d)} \| \sigma_2 \|_{L^2(\mathbb R^n)} \left( \| \Psi_0 \|_{L^2(\mathbb R^n)} + t \|\Psi_1 \|_{L^2(\mathbb R^n)} \right),
\]
for any $0\leq t\leq T$, and, moreover
\[
\| \Phi_0 \|_{C^1([0,T]; L^\infty(\mathbb R^d))} \leq \| \sigma_1 \|_{L^2(\mathbb R^d)} \|\sigma_2\|_{W^{1,2}(\mathbb R^n)} 
\left(2 \| \Psi_0 \|_{L^2(\mathbb R^n)} + (1+T) \|\Psi_1\|_{L^2(\mathbb R^n)} \right). 
\]
\end{itemize}
\end{lemma}

\begin{Proof}
The last statement is a direct consequence of H\"older and Young inequalities; let us detail the proof of items i) and ii).
We associate to  $f\in \big( W^{1,\infty}(\mathbb R^d\times \mathbb R^d) \big)'$, the macroscopic density  $\rho\in \big( W^{1,\infty}(\mathbb R^d) \big)'$ by the formula:
 \[ \left\langle \rho_f , \chi \right\rangle_{(W^{1,\infty})',W^{1,\infty}(\mathbb R^d)} =
   \left\langle f , \chi \otimes \mathbf 1_v \right\rangle_{(W^{1,\infty})',W^{1,\infty}(\mathbb R^d\times \mathbb R^d)}, \qquad \forall \chi \in W^{1,\infty}(\mathbb R^d). \]
Clearly, we have $\| \rho_f \|_{\big( W^{1,\infty}(\mathbb R^d) \big)'} 
 \leq \| f \|_{\big( W^{1,\infty}(\mathbb R^d \times \mathbb R^d) \big)'}$.
 
 For any $\chi \in C^\infty_c(\R^d)$, and  $i\in\{0,1,2\}$ , we can check the following estimates
 \[\begin{array}{lll} \left| \left\langle \rho \ast \Sigma , \nabla^i \chi \right\rangle \right|
 &=& \left| \left\langle \rho , \left( \nabla^i \Sigma \right) \ast \chi \right\rangle \right| 
 \leq \| \rho \|_{\big( W^{1,\infty}(\mathbb R^d) \big)'} \| \left( \nabla^i \Sigma \right) \ast \chi \|_{ W^{1,\infty}(\mathbb R^d)} \\
\qquad\qquad & \leq&  \| f \|_{\big( W^{1,\infty}(\mathbb R^d \times \mathbb R^d) \big)'} 
 \left( \| \nabla^i \Sigma \|_{L^\infty(\mathbb R^d)} +  \| \nabla^{i+1} \Sigma \|_{L^\infty(\mathbb R^d)} \right) \| \chi \|_{L^1(\mathbb R^d)}.  \end{array} \]
Since the dual space of $L^1$ is   $ L^\infty$, for $i=0$, we deduce that
 \[ \begin{array}{lll}
 \| \rho \ast \Sigma \|_{L^\infty(\mathbb R^d)} &\leq& \| f \|_{\big( W^{1,\infty}(\mathbb R^d \times \mathbb R^d) \big)'}  
 \left( \| \Sigma \|_{L^\infty(\mathbb R^d)} +  \| \nabla \Sigma \|_{L^\infty} \right) 
 \\
 &\leq& \| \sigma_1 \|_{W^{1,2}(\mathbb R^d)}^2 \| f \|_{\big( W^{1,\infty}(\mathbb R^d \times \mathbb R^d) \big)'}. \end{array} \]
 Reasoning similarly for $i=1$ and $i=2$, we obtain
 \[ \| \rho \ast \Sigma \|_{W^{2,\infty}(\mathbb R^d)} \leq \|\sigma_1 \|_{W^{3,2}(\mathbb R^d)}^2  \| f \|_{\big( W^{1,\infty}(\mathbb R^d \times \mathbb R^d) \big)'} .\]
We now estimate $p$.
Plancherel's inequality yields
\[ |p'(t)| = \left| \frac1{(2\pi)^n} \int_{\R^n} \cos (c|\xi| t) |\widehat{\sigma_2}(\xi)|^2 \ud \xi \right| \leq \|\sigma_2\|_{L^2(\mathbb R^n)}^2. \]
Since $p(0)=0$, it follows that $|p(t)| \leq  \|\sigma_2\|_{L^2(\mathbb R^n)}^2 t$.
Hence, for all $0 \leq t \leq T<\infty$, we have 
\[ \| \mathcal{L}(f) (t) \|_{W^{2,\infty}(\mathbb R^d\times\mathbb R^d)} \begin{array}[t]{l}
\leq \|\Sigma \ast \rho  \|_{L^\infty(0,T;W^{2,\infty}(\mathbb R^d))} \ds\int_0^t |p(t-s)|  \ud s \\
\leq \|f\|_{C \big( [0,T];\big( W^{1,\infty}(\mathbb R^d\times \mathbb R^d) \big)' \big)} \| \sigma_1 \|^2_{W^{3,2}(\mathbb R^d)} \| \sigma_2 \|^2_{L^2(\mathbb R^n)} \ \ds\frac{T^2}{2}. \end{array}\]
This proves  the estimate in i). That $\mathcal{L}(f)(t)$ is continuous as a function of $t$ follows easily from the previous argument. As a further by-product note that
\[  \| \mathcal{L}(f) (t) \|_{L^\infty} \leq \|f\|_{C \big( [0,T];\big( W^{1,\infty}(\mathbb R^d\times \mathbb R^d) \big)' \big)} \| \sigma_1 \|^2_{W^{1,2}(\mathbb R^d)} \| \sigma_2 \|^2_{L^2(\mathbb R^n)} \ \ds\frac{T^2}{2} \]
holds. Since $p(0)=0$, we have
\[ \partial_t \mathcal{L}(f) (t) =\ds \int_0^t p'(t-s) \Sigma \ast \rho (s) \ud s \]
which gives:
\[  \| \partial_t \mathcal{L}(f) (t) \|_{L^\infty(\mathbb R^d\times \mathbb R^d)} \leq \|f\|_{C \big( [0,T];\big( W^{1,\infty}(\mathbb R^d\times \mathbb R^d) \big)' \big)} \| \sigma_1 \|^2_{W^{1,2}(\mathbb R^d)} \| \sigma_2 \|^2_{L^2(\mathbb R^n)} \ \ds T. \]
This ends the proof of ii).
\end{Proof}

\section{Existence of solutions}\label{Cau}

The proof of existence of solutions to \eqref{simplification} relies on estimates satisfied by  the characteristics curves defined by the following ODE system:
\begin{equation}
\label{char}
\left\lbrace \begin{array}{l} \dot{X}(t) = \xi(t), \\ \dot{\xi}(t) = -\nabla V(X(t)) -\nabla\Phi(t,X(t)). \end{array} \right. 
\end{equation}
From now on, we adopt the following notation.
The potential $\Phi$ being given, we denote by $\varphi_{\alpha}^{\Phi,t}(x_0,v_0)\in \mathbb R^d\times\mathbb R^d$ the solution  of \eqref{char} 
which starts from $(x_0,v_0)$ at time $t=\alpha$:
the initial data is
$\varphi_{\alpha}^{\Phi,\alpha}(x_0,v_0)=(x_0,v_0)$.
We use the shorthand notation $t\mapsto (X(t),\xi(t))$ for $t\mapsto \varphi_{0}^{\Phi,t}(x_0,v_0)$, the solution 
of \eqref{char} with $X(0)=x_0$ and $V(0)=v_0$.
Owing to the regularity of V, $\mathcal L$ and $\Phi_0$, see Lemma~\ref{l:estL}, the 
solution of the differential system \eqref{char}
is indeed well defined for prescribed initial data; this also allows us to establish the following estimates, where characteristics are evaluated both forward and backward.

\begin{lemma}[Estimates on the characteristic curves]\label{l:char}
Let $V$ satisfy \eqref{H2} and let $\Phi\in C^0([0,\infty);W^{2,\infty}(\mathbb R^d))\cap C^1([0,\infty);L^\infty(\mathbb R^d))$.
\begin{enumerate}
\item[a)]
There exists a function $(\mathscr N,t,x,v)\in [0,\infty)\times[0,\infty)\times\mathbb R^d\times \mathbb R^d\mapsto R(\mathscr N,t,x,v)\in [0,\infty)$,
non decreasing with respect to the first two variables, such that the solution $t\mapsto (X(t),\xi(t))$ of \eqref{char} with initial data $X(0)=x_0$, $\xi(0)=v_0$
satisfies the following estimate, for any $t\in \R$,
\[ ( X(t) , \xi(t) ) \in B \left( 0 ,  R\left(\|\Phi\|_{C^1([0,t];L^\infty(\mathbb R^d))}, |t|, x_0, v_0 \right) \right)\subset \mathbb R^d\times\mathbb R^d.  \]
\item [b)]
Taking two different potentials $\Phi_1$ and $\Phi_2$, the following  two estimates hold for any $t>0$:
\[ \begin{array}{l}
 |    (\varphi_0^{\Phi_1,t}-\varphi_0^{\Phi_2,t}) (x_0,v_0)     | 
 \\[.4cm]
 \quad\leq  \ds \int_0^{t} \| (\Phi_1-\Phi_2)(s)\|_{W^{1,\infty}(\mathbb R^d)} 
\exp \left( \int_s^{t} \|  \nabla^2(\Phi_1(\tau)+V)  \| _{L^\infty (B_\tau(x_0,v_0))} \ud \tau \right) \ud  s,
\end{array} \]
\[ \begin{array}{l}
 |    (\varphi_t^{\Phi_1,0}-\varphi_t^{\Phi_2,0}) (x,v)     | 
 \\[.4cm]
 \quad\leq  \ds \int_0^{t} \| (\Phi_1-\Phi_2)(s)\|_{W^{1,\infty}(\mathbb R^d)} 
\exp \left( \int_0^{s} \|  \nabla^2(\Phi_1(\tau)+V)  \| _{L^\infty (\tilde{B}_{t,\tau} (x,v))} \ud \tau \right) \ud  s,
\end{array} \]
where we set $B_\tau(x,v)= B \left( 0 ,  R\left(\max_{i=1,2} \|\Phi_i\|_{C^1([0,\tau];L^\infty(\mathbb R^d))}, \tau , x, v \right) \right)$ and $\tilde{B}_{t,\tau} =  B \left( 0 ,  R\left(\max_{i=1,2} \|\Phi_i\|_{C^1([\tau,t];L^\infty(\mathbb R^d))}, t-\tau , x, v \right) \right)$.

\end{enumerate}
\end{lemma}


%
%

The proof of the lemma is postponed  the end of this section.
Given $0<R_0<\infty$, and $\Psi_0,\Psi_1$ satisfying \eqref{H3} (they enter into the definition of $\Phi_0$ in~\eqref{eq:Phizero}), we set 
\begin{equation} \label{r}
 r(t,x,v)=R(\|\Phi_0 \|_{C^1([0,t];L^\infty(\mathbb R^d))} + |\!|\!| \mathcal{L} |\!|\!|_{\mathcal{B}_t} R_0 ,t,x,v). 
\end{equation}
Proving uniqueness statements for the wide class of external potentials considered in\eqref{H2} requires to strengthen the hypothesis
on the initial data. 
\begin{defin}\label{def1}
Let $0<T,R_0<\infty$.
We say that an integrable function $f_0$ belongs  to the set $E_{R_0,T}$
if  $f_0 \geq 0 $ satisfies $\|f_0\|_{L^1(\mathbb R^d\times\mathbb R^d)}\leq R_0$ and, furthermore,
\[
\mathscr K_{R_0,T}(f_0):=
\ds \int_{\mathbb R^d \times \mathbb R^d} f_0(x,v)
\exp \left( \int_0^T \| \nabla^2 V \|_{L^\infty (B(0 , r(t,x,v)))} \ud t \right) \ud v \ud x<\infty. \]
\end{defin}

\begin{theo}
\label{th1}
Assume \eqref{H1}--\eqref{H3}.
Let $0<R_0,T<\infty$.
Let  $f_0 \in E_{R_0,T}$. Then, there exists a unique $f\in C ([0,T]; L^1 ( \mathbb R^d \times \mathbb R^d ))$ weak solution  of \eqref{simplification}. 
The solution is continuous with respect to the parameters $\mathcal{L}$, $\Phi_0$ and $f_0$, respectively in  $\mathcal{A}_T\cap \mathcal{B}_T$, $C^1([0,\infty);W^{2,\infty}(\mathbb R^d))$ and $E_{R_0,T}$. 
If  $ f_0 \in L^1 (\mathbb R^d \times \mathbb R^d ) $ only, see \eqref{H5}, then there exists  $f\in  C ([0,\infty) ;L^1 ( \mathbb R^d \times \mathbb R^d ))$, weak solution  of \eqref{simplification}.
\end{theo}

\noindent
The statement can be rephrased for the original problem \eqref{kin}--\eqref{CI}.
We also establish the conservation of energy.

\begin{coro}\label{cor1}
Assume \eqref{H1}--\eqref{H3}. Let $0<R_0,T<\infty$.
Let  $f_0 \in E_{R_0,T}$. 
Then, there exists a unique weak solution $(f,\Psi)$ to the system \eqref{kin}--\eqref{CI}
with 
$f\in C([0,T] ; L^1 ( \mathbb R^d \times \mathbb R^d )) $ and $\Psi\in C([0,T]; L^2(\mathbb R^d \times \mathbb R^n))$.
The solution is continuous with respect to the parameters
$\sigma_1$, $\sigma_2$,  $\Psi_0$, $\Psi_1$ and $f_0$
in the sets $W^{3,2}(\mathbb R^d)$,  $L^2(\mathbb R^n)$, $L^2(\mathbb R^d \times \mathbb R^n)$, $L^2(\mathbb R^d \times \mathbb R^n)$ and $E_{R_0,T}$, respectively.
If $ f_0 $ satisfies  \eqref{H5} only, then there exists a weak solution with   $f\in  C ([0,\infty); L^1 ( \mathbb R^d \times \mathbb R^d ))$
and  $\Psi\in C([0,T]; L^2(\mathbb R^d \times \mathbb R^n))$.
Furthermore, when the initial data satisfies  \eqref{H6} the total energy
\[\begin{array}{l}
\ds\frac{1}{2} \int_{\mathbb R^d\times \mathbb R^n} \left| \partial_t \Psi(t,x,y) \right|^2 \ud x  \ud y 
+ \frac{c^2}{2} \int_{\mathbb R^d\times \mathbb R^n} \left| \ds\nabla_y \Psi(t,x,y) \right|^2 \ud x  \ud y 
\\[.4cm]
\qquad\qquad\qquad\qquad\qquad\qquad+ \ds\int_{\mathbb R^d\times \mathbb R^d} f(t,x,v) \left( \frac{|v|^2}{2} + V(x) + \Phi(t,x) \right) \ud v  \ud x
\end{array}\]
is conserved.
\end{coro}

\begin{rmk}
Definition \ref{def1} restricts the set of initial data depending on the growth of the Hessian of the external potential.
Of course, any integrable data $f_0$ with compact support fulfils the criterion in Definition \ref{def1}, 
and when the potential has at most quadratic growth, any data satisfying \eqref{H5} is admissible.
As will be clear in the proof,  the continuity with respect to the initial data does not involve the $L^1$ norm only, but 
 the more intricate quantity $\mathscr K_{R_0,T}$
 also arises in the analysis.
 \end{rmk}

\begin{rmk}
The present approach  does not need  a restriction on the transverse dimension ($n\geq 3$ in \cite{BdB}).
The proof can be slightly modified to treat the case of measure--valued initial data $f_0$, thus including 
the results in \cite{BdB} for a single particle ($f_0(x,v) = \delta_{(x=x_0,v=v_0)}$),
and we can consider a set of $N>1$ particles as well.
The measure--valued solution is then continuous with respect to the initial data in $C([0,T]; (W^{1,\infty}(\mathbb R^d\times\mathbb R^d))')$. 
This viewpoint will be further detailed  with the discussion of  mean--field asymptotics \cite{AVPhD}.
\end{rmk}

The proof of Theorem \ref{th1} relies on a fixed point strategy, the difficulty being to set up the appropriate functional framework.
It turns out that it will be convenient to work with the $C\big([0,T];(W^{1,\infty}(\mathbb R^d \times \mathbb R^d))'\big)$ norm.
We remind the reader that the dual norm on $(W^{1,\infty}(\mathbb R^d \times \mathbb R^d))'$  is equivalent to the Kantorowich--Rubinstein distance
\[W_1(f,g)=\ds\sup_\pi\Big\{\ds\int_{\mathbb R^{2d}\times\mathbb R^{2d}}|\zeta-\zeta'| \ud \pi(\zeta,\zeta')\Big\}\] 
where the supremum is taken over measures $\pi$ having $f$ and $g$ as marginals, see e.~g.~\cite[Remark 6.5]{CV}.
This distance appears naturally in the analysis of Vlasov--like systems, as pointed out in \cite{Dob}.
In order to define the fixed point procedure, we introduce the following mapping.
For a non negative integrable function $f_0$, we denote by  $\Lambda_{f_0}$  the application which associates to 
$\Phi$ in $C([0,\infty);W^{2,\infty}(\mathbb R^d) ) \cap C^1([0,\infty);L^\infty(\mathbb R^d)) $
the unique solution $f$ 
of the Liouville equation
\[
\partial_t f + v\cdot\nabla_x f - \nabla_v f \cdot \nabla_x \left(V + \Phi \right)=0,\]
with initial data $f_0$. We shall make use of the following statement, which provides useful estimates.
\begin{lemma}\label{l:cont}
For any $f_0\in L^1(\mathbb R^d\times \mathbb R^d)$, 
the application  $\Lambda_{f_0}$ is continuous on  
the set $C([0,\infty);W^{2,\infty}(\mathbb R^d) ) \cap C^1([0,\infty);L^\infty(\mathbb R^d)) $ 
with values in 
$C([0,\infty);L^1(\mathbb R^d\times\mathbb R^d)) $.
Furthermore, we have
\[
\|\Lambda_{f_0}(\Phi)-\Lambda_{g_0}(\Phi)\|_{L^\infty(0,\infty;L^1(\mathbb R^d\times \mathbb R^d))}= \|f_0-g_0\|_{L^1(\mathbb R^d\times \mathbb R^d)}
,\]for any $\Phi\in
C([0,\infty);W^{2,\infty}(\mathbb R^d) ) \cap C^1([0,\infty);L^\infty(\mathbb R^d)). $
\end{lemma}

\begin{Proof}
Let $0<T<\infty$ be fixed once for all.
We begin by assuming that  $f_0$ is $C^1$ and compactly supported.
For any $0\leq t\leq T$, we have
\[ \Lambda_{f_0}(\Phi)(t) = f_0 \circ \varphi^{\Phi,0}_{t}, \]
where we remind the reader that $\varphi^{\Phi,0}_{t}(x,v)$ stands for the evaluation at time $0$ of the solution of \eqref{char}
which starts at time $t$ from the state $(x,v)$.
Accordingly any $L^p$ norm is preserved:
$\|\Lambda_{f_0}(\Phi)(t)\|_{L^p(\mathbb R^d\times\mathbb R^d)}=\|f_0\|_{L^p(\mathbb R^d\times\mathbb R^d)}$ holds for any $t\geq 0$ and any $1\leq p\leq \infty$.
By linearity, this immediately proves the continuity estimate with respect to the initial data.

To establish the continuity properties with respect to $\Phi$, we first observe,
denoting  $\Lambda_{f_0}(\Phi)=f$,
 that $(x,v)\in\mathrm{supp}(f(t,\cdot))$ iff $\varphi_t^{\Phi,0}(x,v)\in \mathrm{supp}(f_0)$, that is $(x,v)\in \varphi^{\Phi,t}_0( \mathrm{supp}(f_0))$.
Therefore, by Lemma \ref{l:char},  we can find a compact set $K_T\subset \mathbb R^d\times\mathbb R^d$
such that  $\mathrm{supp}(f(t,\cdot))\subset K_T$ for any $0\leq t\leq T$.
We are dealing with potentials  $\Phi_1$ and $\Phi_2$  in  $C([0,\infty);W^{2,\infty}(\mathbb R^d) ) \cap C^1([0,\infty);L^\infty(\mathbb R^d))$.
We  can again find a compact set, still denoted  by $K_T\subset \mathbb R^d\times\mathbb R^d$, such that the support of the associated solutions 
$ \Lambda_{f_0}(\Phi_1)$ and $\Lambda_{f_0}(\Phi_2)$ 
for any $0\leq t\leq T$ is contained in $K_T$.
We infer that\[\begin{array}{l} \|  \Lambda_{f_0}(\Phi_1)(t)-\Lambda_{f_0}(\Phi_2)(t) \| _{L^1(\mathbb R^d\times\mathbb R^d)}  \ds
= \int_{K_T} 
| f_0 \circ \varphi^{\Phi_1,0}_{t} - f_0 \circ \varphi^{\Phi_2,0}_{t} | \ud v \ud x \\ \ds
\\
\qquad\leq  \|f_0\|_{W^{1,\infty}(\mathbb R^d\times\mathbb R^d)}\ \mathrm{meas}( K_T) \ \ds \sup_{(x,v) \in K_T} |\varphi^{\Phi_1,0}_{t}(x,v) - \varphi^{\Phi_2,0}_{t}(x,v) |
\end{array} \]
holds.
As $\tau$ ranges over $[0,t]\subset[0,T]$ and 
$(x,v)$ lies in $K_T$, 
the backward characteristics $ \varphi^{\Phi_i,\tau }_{t}(x,v)$ still belong to a compact set.
We introduce the following quantities
$$ \mathcal{R} = \sup_{(x,v) \in K_T} R\left(\max_{i=1,2} \|\Phi_i\|_{C^1([0,T];L^\infty(\mathbb R^d))}, T , x, v \right) $$ and
\[ m_T = \exp \left( \int_0^T  \|  \nabla^2\Phi_1(u) \| _{L^\infty( \R^d)}   \ud u \right) .\]
For  $0\leq t\leq T$ and any $(x,v) \in K_T$, Lemma~\ref{l:char}-b) yields:
 \[\begin{array}{l}
 |\varphi^{\Phi_1,0}_{t}(x,v) - \varphi^{\Phi_2,0}_{t}(x,v) | 
 \\[.4cm]
 \qquad\qquad\qquad\leq  m_T\ds \int_0^t \| (\Phi_1-\Phi_2)(s)\|_{  W^{1,\infty}(\mathbb R^d)} 
\exp \left( \int_0^s \|  \nabla^2  V  \| _{L^\infty (B(0,\mathcal{R}))} \ud \tau \right) \ud  s.
\end{array}\]
We conclude with
\[ \sup_{(x,v) \in K_T} |\varphi_t^{\Phi_1,0}(x,v) - \varphi_t^{\Phi_2,0}(x,v) | 
\xrightarrow[\begin{subarray}{c} 
\| \Phi_1-\Phi_2\|_{L^\infty(0,T; W^{2,\infty}(\mathbb R^d))} \rightarrow 0 \\ 
\| \Phi_1 \|_{C^1([0,T];L^\infty(\mathbb R^d))} , \| \Phi_2 \|_{C^1([0,T];L^\infty(\mathbb R^d))} \leq M \end{subarray}
 ]{} 0
. \] 
(It is important to keep both the $C^1([0,T];L^\infty(\mathbb R^d))$ 
and $L^\infty(0,T; W^{2,\infty}(\mathbb R^d))$ norms of the potentials bounded since these quantities appear in the  definition of $\mathcal{R}$ and $m_T$.)
This proves the asserted continuity of the solution with respect to the potential.
By uniform continuity of the flow on the compact set  $[0,T] \times K_T$, we obtain the time continuity.
Hence the result is proved when the initial  data $f_0$ lies  in $C^1_c$.

We finally extend the result for initial data $f_0$ in $L^1$. Those can be approximated by a sequence  $\big(f_0^k\big)_{k\in\mathbb N}$ of functions in $C^1_c(\mathbb R^d \times \mathbb R^d)$.
We have
\[ \| \Lambda_{f_0} (\Phi)(t) - \Lambda_{f_0^k} (\Phi)(t) \|_{L^1(\mathbb R^d\times\mathbb R^d)}  
=  \| \Lambda_{(f_0- f_0^k)} (\Phi)(t) \|_{L^1(\mathbb R^d\times\mathbb R^d)} 
=  \|f_0- f_0^k \|_{L^1(\mathbb R^d\times\mathbb R^d)}. \]
Therefore,
$\Lambda_{f_0}$ is the uniform  limit of maps which are continuous with respect to  $\Phi$ and the time variable.
This remark ends the proof.
\end{Proof}

\begin{ProofOf}{Theorem \ref{th1}}

\noindent
{\it Existence--uniqueness for initial data in $E_{R_0.T}$}. 

We turn to the fixed point reasoning.
For $f$ given in $C\big([0,T];(W^{1,\infty}(\mathbb R^d \times \mathbb R^d))'\big)$, we set
\[ \mathcal{T}_{f_0} (f) = \Lambda_{f_0} (\Phi_0 - \mathcal{L}(f) ). \]
 It is clear that a fixed point of $\mathcal{T}_{f_0}$ is a solution to~\eqref{simplification}. Note also that, as a consequence of Lemma~\ref{l:estL} and Lemma~\ref{l:cont}, $\mathcal{T}_{f_0} (f)(t)\in L^1(\R^d\times\R^d)$. More precisely,  
we know that  $f\mapsto \mathcal{T}(f)$  is continuous with values in  the space $C([0,T]; L^1(\mathbb R^d\times\mathbb R^d)) \subset C\big([0,T];(W^{1,\infty}(\mathbb R^d \times \mathbb R^d))'\big)$.
We shall prove that  $\mathcal{T}$ admits an iteration which is a contraction on the ball with centre  $0$ and radius $R_0$.

Let $f_1$ and $f_2$ be two elements of this ball.
We denote $\varphi_{\alpha}^{\Phi_i,t}$ the flow of \eqref{char} with $\Phi_i = \Phi_0 - \mathcal{L}(f_i)$:
$\varphi^{\Phi_i,t}_{\alpha}(x_0,v_0)$ satisfies \eqref{char} with 
$(x_0,v_0)$ as data at time $t=\alpha$. 
 Let  $ \chi $ be a trial function in  $W^{1,\infty}(\mathbb R^d \times \mathbb R^d)$. We have
\[ \begin{array}{l}
\left| \ds \int_{\mathbb R^d\times\mathbb R^d} (  \mathcal{T}(f_1)(t,x,v) - \mathcal{T}(f_2)(t,x,v) ) \chi(x,v)\ud v\ud x \right| 
\\
\qquad \ds  = \left| \int_{\mathbb R^d\times\mathbb R^d} \left( f_0 \circ \varphi^{\Phi_1,0}_{t} - f_0 \circ \varphi^{\Phi_2,0}_{t} \right)(x,v)
\chi(x,v) \ud v \ud x \right|
\\ \qquad \ds  = \left| \int_{\mathbb R^d\times\mathbb R^d} f_0(x,v) \left( \chi \circ \varphi^{\Phi_1,t}_{0} - \chi \circ \varphi^{\Phi_2,t}_{0} \right)(x,v) \ud v \ud x \right|
\\ \qquad \ds  \leq   \int_{\mathbb R^d\times\mathbb R^d} f_0(x,v)  \|  \nabla \chi \| _\infty 
\left| \varphi^{\Phi_1,t}_{0} - \varphi^{\Phi_2,t}_{0} \right|(x,v) \ud v \ud x. 
\end{array} \]
It follows that
\begin{equation} 
\label{intermediaire}
\|  \mathcal{T}(f_1)(t) - \mathcal{T}(f_2)(t) \| _{\left( W^{1,\infty}(\mathbb R^d\times\mathbb R^d) \right)' }
\leq \int_{\mathbb R^d\times\mathbb R^d} f_0(x,v) 
\left| \varphi^{\Phi_1,t}_{0} - \varphi^{\Phi_2,t}_{0} \right|(x,v) \ud v \ud x .
\end{equation} 
By using  Lemma \ref{l:char}-b), we obtain 
\[\begin{array}{l} \left| \varphi^{\Phi_1,t}_{0} - \varphi^{\Phi_2,t}_{0} \right|(x,v)
\\
 \qquad \leq  \bar m_T\ds \int_0^t \| \mathcal{L}(f_1-f_2)\|_{L^\infty(0,s; W^{2,\infty}(\mathbb R^d))} 
\\
\hspace*{3.5cm}\times
\ds\exp \left( \int_s^t\|  \nabla^2  V  \| _{L^\infty B(0, R(\|\Phi_0 + \mathcal{L}(f_i)\|_{C^1([0,u];L^\infty(\mathbb R^d))} , u, x_0,v_0 ))} \ud u \right) \ud s, \end{array}\]
where we have used
\[\begin{array}{l} \ds \exp \left( \int_0^T  \|  \nabla^2(\Phi_0(u) - \mathcal{L}(f_1)(u) \| _{L^\infty(\mathbb R^d)}   \ud u \right) 
\\[.4cm]
\hspace*{3cm}\leq\ds \exp \left( \int_0^T \left( \|  \nabla^2\Phi_0(u) \| _{L^\infty(\mathbb R^d)} + 
 |\!|\!|  \mathcal{L}  |\!|\!| _{\mathcal{A}_u} \| f_0\| _{L^1(\mathbb R^d\times\mathbb R^d)} \right) \ud u \right)=\bar m_T. \end{array}\] 
Plugging this estimate into \eqref{intermediaire} yields  
\[\begin{array}{l} \|  \mathcal{T}(f_1)(t) - \mathcal{T}(f_2)(t) \| _{\left( W^{1,\infty}(\mathbb R^d\times\mathbb R^d) \right)' }
\\
\qquad\leq \ds \bar m_T  \int_{\mathbb R^d \times \mathbb R^d} f_0(x,v) \int_0^t \|  \mathcal{L}(f_1-f_2)\| _{L^\infty(0,s; W^{2,\infty}(\mathbb R^d))}
\\
\qquad\qquad\qquad\times \ds
\exp \left( \int_s^t \|  \nabla^2 V \| _{L^\infty (B(0 , r(u,x,v)))} \ud u \right)\ud s \ud v \ud x. \end{array} \]
It recasts as 
\[ \|  \mathcal{T}(f_1)(t) - \mathcal{T}(f_2)(t) \| _{\left(W^{1,\infty}\right)'} 
\leq \bar m_T' \mathscr K_{R_0,T}  \int_0^t \| f_1-f_2 \| _{L^\infty \left(0,s; \left(W^{1,\infty}(\mathbb R^d\times \mathbb R^d)\right)'\right)} \ud s \]
with
\[ \bar m_T' = \bar m_T\times    \sup_{0 \leq s \leq T} |\!|\!|  \mathcal{L}  |\!|\!|  _{\mathcal{A}_s} .\]
By  induction, we deduce that
\[ \|  \mathcal{T}^\ell(f_1)(t) - \mathcal{T}^\ell(f_2)(t) \| _{\left(W^{1,\infty}(\mathbb R^d\times \mathbb R^d)\right)'} 
\leq \frac{ \left(t \bar m_T' \mathscr K_{R_0,T} \right)^\ell}{\ell!} \| f_1-f_2 \| _{L^\infty \left(0,T; \left(W^{1,\infty}(\mathbb R^d\times \mathbb R^d)\right)'\right)}  \]
holds  for any $\ell\in\mathbb N$ and $0\leq t\leq T$.
Finally, we are led to
\[ \|  \mathcal{T}^\ell(f_1) - \mathcal{T}^\ell(f_2)\| _{L^\infty \left(0,T; \left(W^{1,\infty}(\mathbb R^d\times \mathbb R^d)\right)'\right)}
\leq \frac{ \left(T \bar m_T' \mathscr K_{R_0,T} \right)^\ell}{\ell!}  \| f_1-f_2 \| _{L^\infty \left(0,T; \left(W^{1,\infty}(\mathbb R^d\times \mathbb R^d)\right)'\right)}.  \]
This shows that  an iteration of   $\mathcal{T}$ is a contraction. Therefore, there exists a unique  fixed point $f$ in $C\big([0,T]; ( W^{1,\infty}(\mathbb R^d\times \mathbb R^d))' \big)$. 
Furthermore, $f=\mathcal T(f)\in
C([0,T]; L^1(\mathbb R^d\times \mathbb R^d))$, and the solution is continuous with respect to the parameters of the system.
Note that the continuity estimate involves the quantity in Definition~\ref{def1} 
which restricts the growth assumption of the initial data.
\medskip

\noindent
{\it Step 2: Existence for an integrable data}

We proceed by approximation. Let $f_0$ be in  $L^1(\mathbb R^d\times\mathbb R^d)$, with $\|f_0\|_{L^1}\leq R_0$.
Then,  
$$(x,v)\mapsto f_0^k(x,v)=f_0(x,v) \mathbf 1_{\sqrt{x^2+v^2}\leq k}$$ lies in $E_{R_0,T}$ (with a constant $\mathscr K_{R_0,T}$ which can blow up as $k\rightarrow\infty$).
The previous step defines $f^k $,  solution of \eqref{simplification} with this initial data. 
Of course we wish to conclude by passing to the limit $k\rightarrow\infty$.
However, the necessary compactness arguments are not direct and the proof splits into several steps.

We start by showing that the sequence $\big(f^k\big)_{k\in\mathbb N}$ is compact  in $C([0,T];\mathcal M^1(\mathbb R^d \times \mathbb R^d)-\text{weak}-\star)$.
Pick $\chi\in C^\infty_c(\mathbb R^d\times\mathbb R^d)$.
For any $0\leq t\leq T$, we have, on the one hand,
\begin{equation}\label{equic}\begin{array}{lll}
 \left| \ds\int_{\mathbb R^d\times\mathbb R^d}f^k (t,x,v) \chi(x,v)\ud v\ud x  \right|
 &\leq &
\| f^k (t,\cdot) \|_{L^1(\mathbb R^d\times\mathbb R^d)} \| \chi \|_{L^\infty(\mathbb R^d\times\mathbb R^d)}
\\
 &\leq & \| f_0^k\|_{L^1(\mathbb R^d\times\mathbb R^d)} \| \chi \|_{L^\infty(\mathbb R^d\times\mathbb R^d)}
\\&\leq &
 \| f_0 \|_{L^1(\mathbb R^d\times\mathbb R^d)} \| \chi \|_{L^\infty(\mathbb R^d\times\mathbb R^d)},\end{array} \end{equation}
and, on the other hand,
\[ \begin{array}{l}
\left| \ds\frac{\ud}{\ud t}\ds\int_{\mathbb R^d\times\mathbb R^d} f^k (t,x,v) \chi(x,v)\ud v\ud x  \right| 
\\
\qquad
= \left| \ds\int_{\mathbb R^d\times\mathbb R^d} f^k (t,x,v) \ \big( v\cdot\nabla_x \chi - \nabla_x ( V+ \Phi_0 - \mathcal{L}(f)(t) )\cdot \nabla_v \chi\big)(x,v) \ud v\ud x \right| \\ \ds 
\qquad\leq  \| f_0 \| _{L^1} 
\Big( \| v\cdot\nabla_x \chi - \nabla V \cdot \nabla_v \chi \|_{L^\infty(\mathbb R^d\times\mathbb R^d))}
\\
\qquad\qquad\qquad
+ \left( |\!|\!| \mathcal{L} |\!|\!|_{\mathcal{A}_T }   \| f_0 \|_{L^1} + \| \Phi_0 \|_{L^\infty([0,T]; W^{1,\infty}(\mathbb R^d))} \Big) 
 \|\nabla_v \chi \|_{L^\infty} \right).
\end{array} \]
Lemma \ref{l:estL} then ensures that the set 
$$\Big\{t\mapsto \ds\ds\int_{\mathbb R^d\times\mathbb R^d} f^k (t,x,v) \chi(x,v)\ud v\ud x,\ k\in\mathbb N\Big\}$$
is equibounded and equicontinuous; hence, by virtue of Arzela--Ascoli's theorem it is relatively compact in $C([0,T])$.
Going back to \eqref{equic}, a simple approximation argument allows us to extend the conclusion to any trial function $\chi$ in  $C_0(\mathbb R^d\times\mathbb R^d)$, the space of continuous functions that vanish at infinity.

This space is separable; consequently, by a diagonal
argument, we can extract a subsequence and find a measure valued function $t\mapsto \ud f(t) \in \mathcal M^1(\mathbb R^d\times\mathbb R^d)$ such that
\[\ds\lim_{k\rightarrow \infty}\ds\int_{\mathbb R^d\times\mathbb R^d} f^k (t,x,v) \chi(x,v)\ud v\ud x
=\ds\int_{\mathbb R^d\times\mathbb R^d}  \chi(x,v) \ud f(t) \]
holds uniformly on $[0,T]$, for any $\chi\in C_0(\mathbb R^d\times\mathbb R^d)$.
As a matter of fact, we note that $\ud f$ is non negative and for any $0\leq t\leq T$ it satisfies
\[\ds\int_{\mathbb R^d\times\mathbb R^d} \ud f(t) \leq \|f_0\|_{L^1(\mathbb R^d\times\mathbb R^d)}. \]

Next, we establish the tightness of the sequence of approximate solutions.
Let  $\epsilon>0$ be fixed once for all.
We can find  $M_\epsilon>0$ such that
\[ \int_{x^2+v^2\geq M_\epsilon^2} f_0(x,v) \ud v \ud x \leq \epsilon. \]
Let us set  $$A_\epsilon= \sup\{ r(T,x,v),\ (x,v) \in B(0,M_\epsilon)\} $$ where we remind the reader that $r(T,x,v)$ has been defined in  \eqref{r}:
 $0<A_\epsilon<\infty$ is well defined by Lemma \ref{l:estL}. 
 Let $\varphi_{\alpha}^{k,t}$ stand for the flow associated to the characteristics of the equation satisfied by $f^k $.
 For any $0\leq t\leq T$, we have
 $\varphi_{0}^{k,t}(B(0,M_\epsilon)) \subset B(0,A_\epsilon)$  so that $\complement \left( \varphi_{t}^{k,0}(B(0,A_\epsilon)) \right)=  \varphi_{t}^{k,0}\big(\complement  B(0,A_\epsilon)\big) \subset \complement B(0,M_\epsilon)$.
It follows that
\[\begin{array}{lll}\ds \int_{\complement B(0,A_\epsilon)} f^k (t,x,v) \ud v \ud x &=& \ds\int_{\complement B(0,A_\epsilon)} f^k _0 ( \varphi_{t}^{k,0}(x,v) ) \ud v \ud x 
\\
&=& \ds\int_{\complement \varphi_{t}^{k,0} (B(0,A_\epsilon))} f_0^k(x,v) \ud v \ud x \\ 
&\leq& \ds\int_{\complement B(0,M_\epsilon)} f_0(x,v) \ud v \ud x \leq \epsilon. 
\end{array} \]
By a standard approximation, we check that the same estimate is satisfied by the limit~$f$:
\[  \int_{\complement B(0,A_\epsilon)} \ud f(t) \leq \epsilon .\]

Finally, we justify that $f^k $ converges to $f$ in $C([0,T]; (W^{1,\infty}(\mathbb R^d \times \mathbb R^d))')$.
Pick $\chi$ in $W^{1,\infty}(\mathbb R^d \times \mathbb R^d)$, with $\|\chi\|_{W^{1,\infty}(\mathbb R^d \times \mathbb R^d)}\leq 1$.
We introduce a cut-off function $\theta_R$ as follows:
\begin{equation}
\label{theta}\begin{array}{ll}
\theta_R(x,v) = \theta(x/R,v/R), \qquad & \theta \in C^\infty_c(\mathbb R^d \times \mathbb R^d),\\
\theta(x,v) =1\text{ for $\sqrt{x^2+v^2}\leq 1$},\qquad &  \theta(x)= 0\text{ for $x^2+v^2\geq 4$}, \\ 0 \leq \theta(x) \leq 1
\text{ for any $x\in\mathbb R^d$}. \end{array}
\end{equation}
Then, we split 
\[ \begin{array}{l}\ds\int_{\mathbb R^d\times\mathbb R^d} f^k (t,x,v)  \chi(x,v) \ud v\ud x -   \ds\int_{\mathbb R^d\times\mathbb R^d}   \chi(x,v) \ud f(t)
\\
\quad=\ds\int_{\mathbb R^d\times\mathbb R^d} f^k (t,x,v)  \chi\theta_R(x,v) \ud v\ud x -   \ds\int_{\mathbb R^d\times\mathbb R^d} \chi\theta_R(x,v) \ud f(t)
\\
\qquad+\ds\int_{\mathbb R^d\times\mathbb R^d} f^k (t,x,v)  \chi(1-\theta_R)(x,v) \ud v\ud x -   \ds\int_{\mathbb R^d\times\mathbb R^d}  \chi(1-\theta_R)(x,v) \ud f(t).
\end{array}\]
Choosing $R\geq A_\epsilon$ yields
\begin{equation}
\label{extraction}\begin{array}{l}
\left| \ds\int_{\mathbb R^d\times\mathbb R^d} f^k (t,x,v)  \chi(1-\theta_R)(x,v) \ud v\ud x -   \ds\int_{\mathbb R^d\times\mathbb R^d}   \chi(1-\theta_R)(x,v) \ud f(t) \right|
\\[.4cm]
\hspace*{10cm}\leq 2 \epsilon \| \chi\| _{L^\infty(\mathbb R^d\times\mathbb R^d)}. 
\end{array}\end{equation}
By virtue of the Arzela-Ascoli theorem, 
$W^{1,\infty}(B(0,2R))$ embeds compactly in  $C (B(0,2R))$. 
Thus, we can find a family $\{\chi_1,..., \chi_{m_\epsilon}\}$ of functions in 
$W^{1,\infty}(\mathbb R^d\times\mathbb R^d)$ such that, for any $\chi\in W^{1,\infty}(\mathbb R^d \times \mathbb R^d)$, $\|\chi\|_{W^{1,\infty}(\mathbb R^d \times \mathbb R^d)}\leq 1$, there exists an index $i\in\{1,...,m_\epsilon\}$ with  $\| \theta_R\chi-\chi_i\| _{L^\infty(B(0,2R))} \leq \epsilon$ (since $\chi\theta_R$ lies in a bounded ball of $W^{1,\infty}(B(0,2R))$). 
Therefore, let us write
\[\begin{array}{l}\ds\int_{\mathbb R^d\times\mathbb R^d} f^k (t,x,v)  \chi\theta_R(x,v) \ud v\ud x -   \ds\int_{\mathbb R^d\times\mathbb R^d}  \chi\theta_R(x,v) \ud f(t)
\\
\qquad= 
\ds\int_{\mathbb R^d\times\mathbb R^d} f^k (t,x,v)  \chi_i(x,v) \ud v\ud x -   \ds\int_{\mathbb R^d\times\mathbb R^d}   \chi_i(x,v) \ud f(t) x
\\
\qquad\qquad+
\ds\int_{\mathbb R^d\times\mathbb R^d} f^k(t,x,v)  (\chi\theta_R-\chi_i)(x,v) \ud v\ud x
-\int_{\mathbb R^d\times\mathbb R^d}  (\chi\theta_R-\chi_i)(x,v) \ud f(t)
,
 \end{array}\]
where the last two terms can both be dominated by $\|f_0\|_{L^1(\mathbb R^d\times\mathbb R^d)}\epsilon$.
We thus arrive at
\[\begin{array}{l}
\left|\ds\int_{\mathbb R^d\times\mathbb R^d} f^k (t,x,v)  \chi(x,v) \ud v\ud x -   \ds\int_{\mathbb R^d\times\mathbb R^d}  \chi(x,v) \ud f(t) \right|
\\
\qquad
\leq 
2 \epsilon (\|\chi\|_{L^\infty(\mathbb R^d\times\mathbb R^d)}+\|f_0\|_{L^1(\mathbb R^d\times\mathbb R^d)})
\\
\qquad\qquad+
\left|\ds\int_{\mathbb R^d\times\mathbb R^d} f^k(t,x,v)\  \chi_i(x,v) \ud v\ud x
-\int_{\mathbb R^d\times\mathbb R^d}  \chi_i(x,v)\ \ud f(t)
\right|
\\
\qquad
\leq 2 \epsilon (\|\chi\|_{L^\infty(\mathbb R^d\times\mathbb R^d)}+\|f_0\|_{L^1(\mathbb R^d\times\mathbb R^d)})
\\
\qquad\qquad+\ds\sup_{j\in\{1,...,m_\epsilon\}}
\left|\ds\int_{\mathbb R^d\times\mathbb R^d} f^k(t,x,v)\  \chi_i(x,v) \ud v\ud x
-\int_{\mathbb R^d\times\mathbb R^d}  \chi_i(x,v)\ \ud f(t)
\right|,
\end{array}
\]
for any $\chi\in W^{1,\infty}(\mathbb R^d \times \mathbb R^d)$, with $\|\chi\|_{W^{1,\infty}(\mathbb R^d \times \mathbb R^d)}\leq 1$.
The last term can be made smaller than $\epsilon$ by choosing $k\geq N_\epsilon$ large enough.
In other words, we can find $N_\epsilon\in\mathbb N$ such that
\[ \begin{array}{l}\ds\sup_{\|  \chi \| _{W^{1,\infty}} \leq 1} \left|\ds\int_{\mathbb R^d\times\mathbb R^d} f^k (t,x,v)  \chi(x,v) \ud v\ud x -   \ds\int_{\mathbb R^d\times\mathbb R^d}  \chi(x,v) \ud f(t)\right|
\\
\qquad\qquad\leq 
2 \epsilon (2+\|f_0\|_{L^1(\mathbb R^d\times\mathbb R^d)})\end{array}\]
holds for any $0\leq t\leq T$, and $k\geq N_\epsilon$:
 $f^k $ converges to  $f$ in
$C\big([0,T]; (W^{1,\infty}(\mathbb R^d \times \mathbb R^d))'\big)$. 
According to Lemma \ref{l:cont}, together with Lemma \ref{l:estL}, it implies that  $\mathcal{T}_{f_0}(f^k )$ converges to $\mathcal{T}_{f_0}(f)$ in 
$C([0,T];L^1(\mathbb R^d\times\mathbb R^d))$. 

By definition  $\mathcal{T}_{f^k _0}(f^k )=f^k $ so that 
\[\begin{array}{l}
 \|  f^k  -\mathcal{T}_{f_0} (f) \| _{C ([0,T];L^1(\mathbb R^d \times \mathbb R^d))}   
\\[.4cm]
\qquad\leq  \|  \mathcal{T}_{f^k _0} (f^k ) -  \mathcal{T}_{f_0} (f^k ) \| _{C([0,T];L^1(\mathbb R^d \times \mathbb R^d)))} 
+\|  \mathcal{T}_{f_0} (f^k ) - \mathcal{T}_{f_0} (f) \| _{C([0,T];L^1(\mathbb R^d \times \mathbb R^d)))}  \\[.4cm] \ds  
\qquad\leq  \|  f^k _0-f_0 \| _{L^1(\mathbb R^d \times \mathbb R^d)} + \|  \mathcal{T}_{f_0} (f^k ) - \mathcal{T}_{f_0} (f) \| _{C([0,T];L^1(\mathbb R^d \times \mathbb R^d)))} 
\end{array} \]
holds, where we have used Lemma \ref{l:cont} again.
Letting $k$ go to $\infty$, we realize that  $f^k $ also converges to $\mathcal{T}_{f_0} (f)$ in  $C([0,T];L^1(\mathbb R^d \times \mathbb R^d))$. 
It implies both $f=\mathcal{T}_{f_0}(f)$ and $f\in C([0,T];L^1(\mathbb R^d \times \mathbb R^d))$. By definition of $\mathcal{T}_{f_0}$, $f$
satisfies  \eqref{simplification}, and it also justifies that $f$ is absolutely continuous with respect to the Lebesgue measure,
which ends the proof. 
\end{ProofOf}

\begin{ProofOf}{Lemma  \ref{l:char}}
Let  $(X,\xi)$ be the solution of  \eqref{char} with $(X(0),\xi(0))=(x_0,v_0)$. We have
\[ \ds\frac{\ud}{\ud t} \left[ V(X(t))+\Phi(t,X(t))+\frac{|\xi(t)|^2}{2} \right] = ( \partial_t \Phi )(t,X(t)). \]
The right hand side is dominated by  $\|  \partial_t \Phi \|_{C([0,t]; L^\infty(\mathbb R^d))} $.
With $t\geq 0$, integrating this relation yields
\[ \ds \frac{|\xi(t)|^2}{2} \leq \Big(V(x_0)+\Phi(0,x_0) + \frac{|v_0|^2}{2} \Big) -(V(X(t))+\Phi(t,X(t))) + t\|  \partial_t \Phi \|_{C([0,t];L^\infty(\mathbb R^d))}. \]
Owing to \eqref{H2} we deduce that 
\[
|\xi(t)|^2
\leq 
a(t) + 2C|X(t)|^2\]
holds with 
\[
a(t)=2 \Big|V(x_0) + \Phi(0,x_0) + \frac{|v_0|^2}{2} \Big|+   2t\|  \partial_t \Phi \|_{C([0,t];L^\infty(\mathbb R^d))}+ 2\| \Phi(t,\cdot) \|_{L^\infty(\mathbb R^d)}+2C.
\]
Next, we simply write
\[
\frac{\ud |X(t)|^2}{\ud t}(t)=2X(t)\cdot \xi(t)\leq X(t)^2+\xi(t)^2\]
so that the   estimate just obtained on $\xi$ yields
\[
|X(t)|^2\leq |x_0| + (1+2C) \ds\int_0^t |X(s)|^2\ud s + \ds\int_0^t a(s) \ud s.
\]
By using the Gr\"onwall lemma we conclude that  
\[ |X(t)|^2 \leq |x_0|^2e^{(1+2C)t} +\ds \int_0^t e^{(1+2C)(t-s)} a(s) \ud s \]
holds.
Going back to the velocity, we obtain
\[ |\xi (t)|^2 \leq 2C \left( |x_0|e^{(1+2C)t} + \ds\int_0^t  e^{(1+2C)(t-s)} a(s) \ud s  \right) + a(t). \]
It concludes the proof of Lemma~\ref{l:char}-a).
\\

Next, let  $(X_1,\xi_1)$ and $(X_2,\xi_2)$ be two solutions of \eqref{char} with the same initial data $(x_0,v_0)$, but different  potentials 
$\Phi_1,\Phi_2$. We already know that the two characteristic curves $(X_i(s),\xi_i(s))$, for $i\in\{1,2\}$,
 belong to $B_s(x,v)$. 
We have
\[ \left\lbrace \begin{array}{l}
\ds\frac{\ud}{\ud s}|X_1(s)-X_2(s)| \leq |\xi_1(s)-\xi_2(s)|, \\[.4cm]
\ds\frac{\ud}{\ud s}|\xi_1(s)-\xi_2(s)| \leq
\| \nabla \left( \Phi_1(s,\cdot) - \Phi_2(s,\cdot) \right) \|_{L^\infty(\mathbb R^d)} \\
\hspace*{4cm}+ |X_1(s)-X_2(s)| 
\|  \nabla^2 ( V+\Phi_1(s,\cdot) ) \| _{L^\infty (B_s(x,v))} \end{array} \right. \] 
The Gr\"onwall lemma yields the estimate 
\[\begin{array}{l} |(X_1(t),\xi_1(t))-(X_2(t),\xi_2(t))| 
\\
\leq\ds  \int_0^t  \|(\Phi_1-\Phi_2)(\tau,\cdot)\|_{W^{1,\infty}(\mathbb R^d)}  \exp \left( \int_s^t
\big(\|  \nabla^2(V+\Phi_1(u)) \| _{L^\infty (B_u(x,v))}\big) \ud u \right) \ud s. \end{array}\]
Finally, we wish to evaluate the backward characteristics, looking at the state at time 0, 
given the position/velocity pair at time $t$.  
Namely we consider $\varphi_t^{\Phi,s}(x,v)$ for $s
\leq t$, bearing in mind $ \varphi_t^{\Phi,t}(x,v)=(x,v)$.
We set 
 $$\begin{pmatrix}Y \\ \zeta\end{pmatrix}(s) =\begin{pmatrix} 1 & 0
 \\ 
 0 & -1\end{pmatrix} \varphi_t^{\Phi,t-s}(x,v).$$ We  check that $(Y,\zeta)$ satisfies
\[\left\lbrace \begin{array}{ll} \ds\frac{\ud}{\ud s}Y(s) = \zeta(s), \qquad &  \ds\frac{\ud}{\ud s}\zeta(s) = -\nabla  V(Y(s)) -\nabla\Phi(t-s,Y(s)), 
\\ [.4cm]
Y(0) = x, \qquad & \zeta(0) = v. \end{array} \right.
 \]
Changing $\Phi$ for $\Phi(t-\cdot)$, this allows
 us to obtain  the same estimates on $(Y,\zeta)$ for all $s\geq0$. We conclude by taking $s=t$.
\end{ProofOf}


\begin{ProofOf}{Corollary \ref{cor1}}
Theorem \ref{th1}  constructs solutions to \eqref{simplification} in $C^0([0,\infty);L^1(\mathbb R^d\times\mathbb R^d))$. We have now the functional framework
 necessary to justify the manipulations made in Section~\ref{s:prelim}.
For $\Psi_0,\Psi_1$ verifying \eqref{H3}, 
formula \eqref{Fourrier} defines a solution $\Psi\in C([0,\infty);L^2(\mathbb R^n\times\mathbb R^d))$ of the wave equation, and finally $(f,\Psi)$ satisfies \eqref{kin}--\eqref{CI}.
Conversely, if $ f\in C^0([0,\infty);L^1(\mathbb R^d\times\mathbb R^d))$ and $\Psi\in C([0,\infty);L^2(\mathbb R^n\times\mathbb R^d))$ is a solution of the system \eqref{kin}--\eqref{CI},
then we can rewrite $\Phi=\Phi_0-\mathcal L(f)$ and $f$ verifies \eqref{simplification}.
This equivalence justifies the first part of the statement in Corollary~\ref{cor1}.

It only remains to justify the energy conservation.
We consider an initial data with finite energy: 
\[\begin{array}{l}\mathcal{E}_0 =\underbrace{\ds \frac{c^2}{2} \int_{\mathbb R^d\times \mathbb R^n} \left| \nabla_y \Psi_0 (x,y) \right|^2 \ud y \ud x 
+ \ds \frac{1}{2} \int_{\mathbb R^d\times \mathbb R^n} \left| \Psi_1(x,y) \right|^2 \ud y \ud x }_{\mathcal E_{0}^{\mathrm{vib}}}
\\
\qquad\qquad+\underbrace{\ds \int_{\mathbb R^d\times \mathbb R^d} f_0(x,v) \left( \frac{|v|^2}{2} + V(x) + \Phi(0,x) \right)\ud v \ud x}_{\mathcal E_0^{\mathrm{part}}} \in (-\infty,+\infty).\end{array}\]
For the solutions constructed in Theorem \ref{th1}, we have seen that 
the self--consistent potential remains smooth enough so that the characteristic curves $t\mapsto (X(t),\xi(t))$ are well--defined.
 Therefore, we can write
 \[\begin{array}{l} \ds\int_{\mathbb R^d\times \mathbb R^d} f(t,x,v) \left( \frac{|v|^2}{2} + V(x) + \Phi(t,x) \right) \ud v \ud x
\\
\qquad\qquad= \ds\int_{\mathbb R^d\times \mathbb R^d} f_0(x,v) \left( \frac{|\xi(t)|^2}{2} + V(X(t)) + \Phi(t,X(t)) \right) \ud v \ud x.\end{array}\]
For any $(t,x,v)$ we have the following equality
\[ \ds\frac{\ud}{\ud t} \left[ V(X(t))+\Phi(t,X(t))+\frac{|\xi(t)|^2}{2}  \right] = ( \partial_t \Phi )(t,X(t)). \]
Therefore, we get
\[ \begin{array}{l} \ds\int_{\mathbb R^d\times \mathbb R^d} f(t, x, v) \left( \frac{|v|^2}{2} + V(x) + \Phi(t,x) \right) \ud v \ud x \ds 
\\[.4cm]
\qquad\qquad=\ds \mathcal{E}_0^{\mathrm{part}}
+ \int_{\mathbb R^d\times \mathbb R^d} f_0(x,v) \int_0^t ( \partial_t \Phi )(s,X(s)) \ud s \ud v \ud x
\\ [.4cm]\ds 
\qquad\qquad= \mathcal{E}_0^{\mathrm{part}}
+ \int_0^t \int_{\mathbb R^d\times \mathbb R^d} f(s,x,v) (\partial_t \Phi )(s,x) \ud v \ud x \ud s \\ [.4cm]\ds 
\qquad\qquad= \mathcal{E}_0^{\mathrm{part}}
+ \int_0^t \int_{\mathbb R^d} \rho(s,x) (\partial_t \Phi )(s,x) \ud x \ud s. \end{array} \]
Next, let $\Psi$ be the unique solution of  \eqref{modele_milieu} associated to $f$. We first assume that the initial data  $\Psi_0$ et $\Psi_1$ are smooth, say in  $L^2(\mathbb R^d,H^2(\mathbb R^n))$.
Therefore, going back to 
  \eqref{Fourrier}, we can check that $\Psi$ lies in 
$C([0,\infty);L^2(\mathbb R^d,H^2(\mathbb R^n)))$. Integrations by parts lead to
\[ \begin{array}{l}\ds\frac{\ud}{\ud t} \left[  \frac{1}{2} \int_{\mathbb R^d\times \mathbb R^n} \left| \partial_t \Psi(t,x,y) \right|^2\ud y  \ud x 
+ \frac{c^2}{2} \int_{\mathbb R^d\times \mathbb R^n} \left| \nabla_y \Psi(t,x,y)  \right|^2 \ud x \ud y \right]  \ds 
\\[.4cm]
\qquad\qquad\ds
=  \int_{\mathbb R^d\times \mathbb R^n} \partial_t \Psi \left( \partial_t^2 \Psi - c^2 \Delta_y \Psi \right)t,x,y)\ud y  \ud x \\ [.4cm]\ds 
\qquad\qquad
= -\int_{\mathbb R^d\times \mathbb R^n} \partial_t \Psi(t,x,y)\ \rho(t,\cdot) \underset{x}{\ast} \sigma_1(x)\  \sigma_2 (y) \ud y \ud x   \\ [.4cm]\ds 
\qquad\qquad= -\int_{\mathbb R^d} \rho \partial_t \Phi(t,x) \ud x. \end{array}\]
Hence, we obtain
\[\begin{array}{l}\ds\frac{1}{2} \int_{\mathbb R^d\times \mathbb R^n} \left| \partial_t \Psi (t,x,y)\right|^2 \ud x \ud y 
+ \frac{c^2}{2} \int_{\mathbb R^d\times \mathbb R^n} \left| \nabla_y \Psi (t,x,y) \right|^2 \ud x \ud y
\\[.4cm]
\qquad\qquad\ds=  \mathcal{E}_0^{\mathrm {vib}} - \int_0^t \int_{\mathbb R^d} \rho(s,x) (\partial_t \Phi )(s,x) \ud x \ud s.\end{array} \]
It proves the energy conservation for such smooth data.

We go back to general data with finite energy: $\Psi_0\in L^2(\mathbb R^d,H^1(\mathbb R^n))$ and $\Psi_1\in L^2(\mathbb R^d \times \mathbb R^n)$. We approximate the data  by $\Psi_0^k$ and 
$\Psi_1^k$ lying  in  $L^2(\mathbb R^d,H^2(\mathbb R^n))$. 
Using  \eqref{Fourrier}, one sees the associated sequence $(\Psi^k)_{k\in\mathbb N}$ of solutions to  \eqref{modele_milieu} converges to $\Psi$ in  $ C([0,\infty);L^2(\mathbb R^d, H^1(\mathbb R^n)))$ and 
$ C^1([0,\infty);L^2(\mathbb R^d \times \mathbb R^n))$. This implies one can pass to the limit in the energy conservation.
\end{ProofOf}

\begin{rmk}
We point out that, whereas energy conservation 
is an important physical property, it was not used here in the existence proof.
In particular, one should notice that it does not provide directly useful a 
priori estimates on the kinetic energy, since the potential energy associated to the external potential $V$ can be negative and unbounded under our assumptions.
In order to deduce a useful estimate the assumptions on the initial data need to be strengthened:
in addition to \eqref{H6} we suppose
\[M_2:=\ds \int_{\mathbb R^d\times \mathbb R^d} f_0(x,v) |x|^2 \ud v \ud x< \infty. \]
We set $V_-(x)=\max (-V(x),0)\geq 0$.
Then \eqref{H2} implies
\[ \int_{\mathbb R^d\times \mathbb R^d} f(t,x,v) V_-(x) \ud v \ud x\begin{array}[t]{l} \ds \leq  \int_{\mathbb R^d\times \mathbb R^d} f(t,x,v) C(1+|x|^2) \ud v \ud x\\ \ds 
\leq  C\|f_0\|_{L^1(\mathbb R^d\times \mathbb R^d)} + C \int_{\mathbb R^d\times \mathbb R^d} f_0(x,v) 
|X(t)|^2 \ud v \ud x,\end{array} \]
where  $X(t)$ stand for the first (space)  component of $\varphi_0^t(x,v)$.
Reproducing the estimates of the proof of Lemma \ref{l:char}, we get 
\[ 
|X(t)|
\leq |x|e^{\sqrt{2 C}t} + \frac{1}{\sqrt{C}} \left( V(x) + \frac{|v|^2}{2} + \Phi(0,x) \right)^{1/2} (e^{\sqrt{2 C}t}-1) + b(t) \]
where
\[ b(t) = \sqrt{2} \int_0^t \left( C+ \|\Phi(s,\cdot)|_{L^\infty(\mathbb R^d)} + s \| \partial_t \Phi \|_{ C([0,s];L^\infty(\mathbb R^d))} \right)^{1/2} 
e^{\sqrt{2 C}(t-s)} \ud s .\]
It follows that 
\[ 
|X(t)|
\leq 9 |x|^2 e^{2\sqrt{2 C}t} + \frac{9}{C} \left( V(x) + \frac{|v|^2}{2} + \Phi(0,x) \right) (e^{ \sqrt{2 C}t}-1)^2 + 9 b(t)^2. \]
Eventually, we find
\[
\int_{\mathbb R^d \times \mathbb R^d} f(t,x,v) V_-(x) \ud v \ud x
\leq
  Ce^{2\sqrt{2 C}t}M_2
+ 9(e^{ \sqrt{2 C}t}-1)^2 \mathcal{E}_0 + C(9b(t)^2+1) \|f_0\|_{L^1(\mathbb R^d\times \mathbb R^d)}. \]
Therefore the potential energy associated to the external potential cannot be too negative and all terms in the energy balance remain bounded on any finite time interval.
\end{rmk}

\section{Large wave speed asymptotics}\label{As}

This section is devoted to the asymptotics of large wave speeds.
Namely, we consider the following rescaled version of the system:
\begin{equation}
\label{limite1}
\left\lbrace \begin{array}{l}
\partial_t f_\epsilon + v\cdot\nabla_x f_\epsilon -\nabla_x (V+\Phi_\epsilon)\cdot \nabla_v f_\epsilon =0, \\
\Phi_\epsilon(t,x,y)=\ds\int_{\mathbb R^n\times\mathbb R^d} \Psi_\epsilon(t,z,y)\sigma_2(y)
\sigma_1(x-z)\ud z\ud y,
\\
\Big(\partial^2_{tt}-\ds\frac1\epsilon\Delta_y \Big)\Psi_\epsilon(t,x,y)=-\ds\frac1\epsilon
\sigma_2(y)\ds\int_{\mathbb R^d\times\mathbb R^d} \sigma_1(x-z)f(t,z,v)\ud v\ud z,
\end{array} \right.
\end{equation}
completed with suitable initial conditions. We are interested in the behavior of the solutions as $\epsilon\rightarrow 0$.
We shall discuss below the physical meaning of this regime.
But, let us first explain on formal grounds what can be expected.
As $\epsilon\rightarrow 0$ the wave equation degenerates to
\[- \Delta_y \Psi(t,x,y)=-
\sigma_2(y)\ \sigma_1\underset{x}{\ast}\rho(t,x),\qquad \rho(t,x)=\ds\int_{\mathbb R^d}f(t,x,v)\ud v.\]
We obtain readily 
the solution by uncoupling the variables:
$$\Psi(t,x,y)=\gamma(y)\ \sigma_1\underset{x}{\ast}\rho(t,x)$$ where $\gamma$ satisfies the mere Poisson equation
$\Delta_y \gamma=\sigma_2$.
At leading order the potential then becomes
\[
\Phi(t,x)=-\kappa \ \Sigma\underset{x}{\ast}\rho(t,x),\qquad\Sigma=\sigma_1\ast\sigma_1,
\qquad
\kappa =- \ds\int_{\mathbb R^n} \sigma_2\gamma\ud y.\]
Therefore, we guess that the limiting behavior is described by the following Vlasov equation
\[\partial_t f + v\cdot\nabla_x f -\nabla_x (V+\Phi)\cdot \nabla_v f =0.\]
As long as the integration by parts makes sense (we shall see that difficulties in the analysis precisely arise when $n\leq 2$), we observe that
\[\kappa = \ds\int_{\mathbb R^n} |\nabla_y \gamma|^2\ud y>0. \]
It is then tempting to make the form function   $\sigma_1$ 
depend on $\epsilon$ too, so that $\Sigma$ resembles the kernel of $(-\Delta _x)$.
We would arrive at the Vlasov--Poisson system, in the case of attractive forces.
We which to justify such asymptotic behavior.

\subsection{Dimensional analysis}

In \eqref{kin}, $f$ is the density of particles in phase space: it gives a number of particles per unit volume
of phase space.
Let $T,L,\mathcal V$ be units for time, space and velocity respectively, and set
\[t'=t/T,\quad x'=x/L,\quad v'=v/\mathcal V\]
which define dimensionless quantities.
Then, we set
\[f'(t',x',v')\ L^{-d}\ \mathcal V^{-d}=f(t,x,v)\]
(or maybe more conveniently $f'(t',x',v')\ud v'\ud x'=f(t,x,v)\ud v\ud x$).
The external and interaction potential, $V$ and $\Phi$, have both the dimension of a velocity squared.
We set
\[V(x)=\mathcal V_{\mathrm{ext}}^2\ V'(x'),\qquad 
\Phi(t,x)=\mathcal W^2\ \Phi'(t',x'),\]
where $\mathcal V_{\mathrm{ext}}$ and $\mathcal W$ thus have the dimension of a velocity.
We switch to the dimensionless equation
\[
\partial_{t'} f'+\ds\frac{\mathcal V T}{L}v'\cdot\nabla_{x'} f'-\ds\frac{ T}{L\mathcal V}\mathcal V^2\nabla_{x'}\Big(V'+\Big(\ds\frac{\mathcal W}{\mathcal V}\Big)^2\Phi'\Big)\cdot\nabla_{v'} f'=0.
\]
The definition of the interaction potential $\Phi$ is driven by the product $\sigma_2(z)\sigma_1(x)\ud x$.
We scale it as follows
\[\sigma_2(z)\sigma_1(x)\ud x=\Sigma_\star L^d \sigma_2'(z')\sigma_1'(x')\ud x'.\]
It might help the intuition to think $z$ as a length variable, and thus $c$ has a velocity, but there is not reason to assume such privileged units.
Thus, we keep a general approach.
 For the vibrating field, we set
 \[
 \psi(t,x,z)=\Psi_\star \ \psi'(t',x',z'),\qquad z'=z/\ell,
 \]
 still with the convention that primed quantities are dimensionless.
 Accordingly, we obtain
 \[
 \mathcal W^2=\Sigma_\star L^d\Psi_\star\ell^n\]
 and the consistent expression of the dimensionless potential
 \[
 \Phi'(t',x')=\ds\int  \sigma_1'(x'-y')\sigma_2'(z')\psi(t',y',z')\ud z'\ud y'.\]
 The wave equation becomes
 \begin{equation}\label{eq:scalewave}
 \partial^2_{t't'}\psi'-\ds\frac{T^2 c^2}{\ell^2}\Delta_{z' }\psi'=-\underbrace{\ds\frac{T^2 \Sigma_\star L^d}{\Psi_\star} L^{-d}}_{\ds\frac{T^2 \Sigma_\star}{\Psi_\star}}\ \sigma_2'(z')\ds\int \sigma_1'(x'-y')f'(t',y',v')\ud v'\ud y'.
 \end{equation}
 Note that 
 \[\ds\frac{T^2 \Sigma_\star}{\Psi_\star}=\Sigma_\star L^d\ell^n\Psi_\star\ \ds\frac{T^2}{\Psi_\star^2L^d\ell^n}=\mathcal W^2 \ \ds\frac{T^2}{\Psi_\star^2L^d\ell^n}.\]
 
Let us consider the energy 
 balance where the following quantities, all having the homogeneity of a velocity squared,  appear:
 \begin{itemize}\item
 the kinetic energy of the particles $\int v^2 f\ud v\ud x$; it scales like $\mathcal V^2$,
 \item the external potential energy $\int Vf\ud v\ud x$; it scales like $\mathcal V_{\mathrm{ext}}^2$,
 \item the coupling energy $\int \Phi f\ud v\ud x$; it scales like $\mathcal W^2$,
 \item the wave energy which splits into:
  \begin{itemize}\item[a)]
  $\int |\partial_t\psi|^2\ud z\ud x$, which scales like $\Psi_\star^2\frac{L^d\ell^n}{T^2}$,
   \item[b)] $c^2\int |\nabla_z\partial_t\psi|^2\ud z\ud x$, which scales like $c^2\Psi_\star^2\frac{L^d\ell^n}{\ell^2}$.
 \end{itemize}
 Note that the kinetic energy in a) is $\frac{\ell^2 }{c^2T^2}$ times the elastic energy in b).
  \end{itemize}
To recap, we have at hand 
6 parameters imposed by the 
model ($L,  \ell,  c, \mathcal V_{\mathrm{ext}}, \mathcal W, \Sigma_\star$)
and two parameters governed by the initial conditions $\mathcal V$ and $\Psi_\star$. 
They allow to define the five energies described above.
\\

We turn to the scaling assumptions.
It is convenient to think of them by comparing the different time scales involved in the equations.
We set 
\[\epsilon=\left(\ds\frac{\ell}{c T}\right)^2\ll 1.\]
If $\ell$ is the size of the support of the source $\sigma_2$, 
then this regime means that the time a typical
particle needs to cross $ L$ (the support of $\sigma_1$) is much longer than the time the wave
needs to cross $\ell$ (the support of $\sigma_2$). 
Next we suppose that the kinetic energy of the particle, the 
energy of the particle associated to the external potential,  the elastic energy of
the wave as well as the interaction energy, all have the same strength, which can expressed by setting
\[\ds\frac{L}{T}= \mathcal V=  \mathcal V_{\mathrm{ext}}= \mathcal W= \sqrt{c^2\Psi_\star^2 L^d\ell^{n-2}}.
\]
As a consequence, it imposes the following scaling of the coupling constant
\[\ds\frac{\Psi_\star}{T^2\Sigma_\star}=\epsilon.
\]
It also means that the kinetic energy of the wave is small 
  with respect to its elastic energy. Inserting this in~\eqref{eq:scalewave} yields~\eqref{limite1}.

\subsection{Statements of the results}

Throughout this Section, we assume \eqref{H1}, and 
we shall strengthen the assumptions {\eqref{H2}--\eqref{H6} as follows (note that since we are dealing with sequences of initial data, it is important to make
 the estimates uniform with respect to the scaling parameter):


\begin{equation}
\tag{{\bf{H7}}}\label{H7}
\textrm{the external potential $V\in W^{2,\infty}_{\mathrm{loc}}(\mathbb R^d)$ is non negative,}
\end{equation}
\begin{equation}
\tag{{\bf{H8}}}\label{H8}
\left.\begin{array}{l}
\textrm{$f_{0,\epsilon}\in L^1(\mathbb R^d\times\mathbb R^d)$, with a uniformly bounded norm, 
}
\\
\textrm{and $\Psi_{0,\epsilon},\Psi_{1,\epsilon}\in L^2(\mathbb R^d\times\mathbb R^n)$
are such that the rescaled initial energy}
\\
\mathcal E_{0,\epsilon}
=\ds\int_{\mathbb R^d\times\mathbb R^d}\Big(\ds\frac{v^2}{2} +V+|\Phi_\epsilon| \Big)f_{0,\epsilon}\ud v\ud x
\\
\qquad\qquad\qquad
+\ds\frac\epsilon2\ds\int_{\mathbb R^n\times\mathbb R^d}| \Psi_{1,\epsilon}|^2\ud y\ud x 
+\ds\frac1{2}\ds\int_{\mathbb R^n\times\mathbb R^d}|\nabla_y \Psi_{0,\epsilon}|^2\ud y\ud x 
\\
\textrm{is uniformly bounded: $0\leq \sup_{\epsilon>0} \mathcal E_{0,\epsilon}=\bar{\mathcal E}_0<\infty$.}
\end{array}\right\}\end{equation}
\begin{equation}
\tag{{\bf{H9}}}\label{H9}
\textrm{$f_{0,\epsilon}$ is bounded in $ L^\infty(\mathbb R^d\times\mathbb R^d)$, uniformly with respect to $\epsilon$.} 
\end{equation}

\begin{theo}\label{cvth1}
Suppose $n \geq 3$. Let \eqref{H1} and \eqref{H7}--\eqref{H9} be satisfied.
Let $(f_\epsilon,\Psi_\epsilon)$ be the associated solution to \eqref{limite1}.
Then, there exists a subsequence such that $f_\epsilon$ converges in $C([0,T];
L^p((\mathbb R^d\times\mathbb R^d)-\text{weak}))$ for any $1\leq p<\infty$  
 to $f$ solution of the following Vlasov equation  
\begin{equation}
\label{Vlasov}
\left\lbrace \begin{array}{l}
\partial_t f + v\cdot\nabla_x f - \nabla_x (V +\bar\Phi)\cdot\nabla_v f  =0, \\
f(0,x,v)= f_0(x,v), \end{array} \right.
\end{equation}
where
\[ 
\bar\Phi=- \kappa \Sigma \ast \rho ,\qquad \Sigma=\sigma_1\underset{x}{\ast}\sigma_1,\qquad
\kappa = \int_{\mathbb R^n} \frac{|\widehat{\sigma_2}(\xi)|^2}{(2\pi)^n|\xi|^2} \ud \xi, \]
and $f_0$ is the weak 
 limit in $L^p(\mathbb R^d\times\mathbb R^d)$ 
of $f_{0,\epsilon}$. 
\end{theo}

In order to derive the Vlasov--Poisson system from \eqref{limite1}, the form function $\sigma_1$ need to be appropriately defined and scaled with respect to $\epsilon$.
Let  $\theta$ and $\delta$ be two radially symmetric functions in  $C^\infty_c(\mathbb R^d) $ verifying:
\[ 0 \leq \theta , \delta \leq 1 \qquad  \theta(x) = 1 \mbox{ for $|x|\leq 1$}, \qquad  \theta (x) =0 \mbox{ for $|x|\geq 2$},
\qquad  \int_{\mathbb R^d} \delta(x) \ud x =1. \]
We set $\theta_\epsilon(x) = \theta(\sqrt\epsilon x)$ et $\delta_\epsilon(x)= \ds\frac{1}{\epsilon^{d/2}} \delta( x/\sqrt\epsilon)$ and
\[ \sigma_{1,\epsilon} = C_d \delta_\epsilon \ast \frac{\theta_\epsilon}{|\cdot|^{d-1}}, \qquad
\textrm{with } 
C_d=\left(|\mathbb S^{d-1}|\ds\int_{\mathbb R^d}\ds\frac{\ud x}{|x|^{d-1}|e_1-x|^{d-1}}\right)^{-1/2}.\]

\begin{theo}\label{cvth2}
Let  $d=3$ and $n \geq 3$. Assume \eqref{H1} and \eqref{H7}--\eqref{H9}. 
Let $(f_\epsilon,\Psi_\epsilon)$ be the associated solution to \eqref{limite1}.
Then, there exists a subsequence such that $f_\epsilon$ converges in $C([0,T];L^p(\mathbb R^3\times\mathbb R^3)-\text{weak})$ for any $1< p<\infty$
 to $f$ solution of the attractive  Vlasov--Poisson equation   
 \begin{equation}
\label{VlasovP}
\left\lbrace \begin{array}{l}
\partial_t f + v\cdot\nabla_x f - \nabla_x (V +\bar\Phi)\cdot\nabla_v f  =0, \\
\Delta\bar \Phi=\kappa \rho,
\\
f(0,x,v)= f_0(x,v) \end{array} \right.
\end{equation}
where $f_0$ is the weak limit in $L^p(\mathbb R^3\times\mathbb R^3)$ of $f_{0,\epsilon}$.
\end{theo}

\begin{rmk}\label{plus:cv}
In Theorem \ref{cvth1}, if, furthermore, we assume that 
$\big(f_{0,\epsilon}\big)_{\epsilon>0}$ 
converge
 (in the appropriate weak sense)  to $f_0$, by uniqueness of the solution of the limit equation, the entire sequence 
$\big(f_\epsilon\big)_{\epsilon>0}$ converges to $f$.
For Theorem \ref{cvth1} and Theorem \ref{cvth2}, if the initial data converges strongly to $f_0$ in $L^p(\mathbb R^d\times\mathbb R^d)$, $1\leq p<\infty$,
then $f_\epsilon$ converges to $f$ in $C([0,T];L^p( \mathbb R^d\times\mathbb R^d))$.
\end{rmk}

\subsection{Convergence to the Vlasov equation with a smooth convolution kernel}

Taking into account the rescaling, the analog of \eqref{simplification} for \eqref{limite1} reads
\begin{equation}\label{simpl2}
\partial_t f_\epsilon+ v\cdot\nabla_x f_\epsilon- \nabla_x\Big(V+\Phi_{0,\epsilon}-\ds\frac1\epsilon \mathcal L_\epsilon(f_\epsilon)\Big)\cdot\nabla_vf_\epsilon=0,
\end{equation}
with
 \[ \Phi_{0,\epsilon}(t,x) = \int_{\R^d \times \R^n} \widetilde \Psi_{\epsilon}(t,z,y) \sigma_1(x-z)  \sigma_2(y) \ud y \ud z. \] 
where  $\widetilde\Psi_{\epsilon}$ stands for  the unique solution of the free linear wave equation (in $\mathbb R^n$) with
wave speed $1/\epsilon$ and  initial data $\Psi_{0,\epsilon}$ and $\Psi_{1,\epsilon}$, 
and
\begin{equation}\label{resc_pot}\begin{array}{lll}
\ds\frac1\epsilon \mathcal L_\epsilon(f_\epsilon)(t,x)&=&\ds\frac1\epsilon \ds\int_{\mathbb R^d}\Sigma(x-z)\left(  \int_0^t \rho_\epsilon(t-s,z) 
\right.\\
[.4cm]
&&\left.\qquad
\qquad\qquad\times\left(
\ds\int_{\mathbb R^n} \ds\frac{\sin(|\xi|s/\sqrt\epsilon)}{|\xi|/\sqrt\epsilon}\  |\widehat{\sigma}_2(\xi)|^2 \ds\frac{\ud \xi}{(2\pi)^n} \right) \ud s\right)\ud z
 \\[.4cm]
 & =&\ds\left( \Sigma
 \underset{x}{\ast}  \int_0^{t/\sqrt\epsilon} \rho_\epsilon(t-s\sqrt\epsilon,\cdot) \ q(s) \ud s\right)(x)
\end{array}\end{equation}
where we have set
\[q(t)=\ds\frac{1}{(2\pi)^n} \int_{\mathbb R^n} \frac{\sin(t |\xi|) }{|\xi|} |\widehat{\sigma_2}(\xi)|^2 \ud\xi.\]
(it is nothing but $p(t)$ as introduced in Section~\ref{s:prelim} evaluated with $c=1$; of course when $c=1$ and $\epsilon=1$, the operators $\frac1\epsilon\mathcal L_\epsilon$ 
in \eqref{resc_pot} and $\mathcal L$ in \eqref{eq:Loperator} coincide.)

\begin{lemma}\label{qL1}
Let $n\geq 3$. Then $q$ is integrable over $[0,+\infty[$ with 
\[ \int_0^\infty q(t) \ud t =\ds\frac{1}{(2\pi)^n} \int_{\mathbb R^n} \frac{|\widehat{\sigma_2}(\xi)|^2 }{|\xi|^2}  \ud\xi:= \kappa >0. \]
\end{lemma}

\begin{Proof}
By virtue of the dominated convergence theorem, 
$t\mapsto q(t)$ is  continuous on $[0,\infty)$. Bearing  in mind that $\sigma_2$ is radially symmetric, integrations by parts yield
\[ \begin{array}{lll}q(t) &=& \ds \frac{|\mathbb{S}^{n-1}|}{(2\pi)^n} \int_0^\infty \sin(tr) r^{n-2} |\widehat{\sigma_2}(re_1)|^2 \ud r 
\\[.4cm] \ds 
&=&\ds \frac{|\mathbb{S}^{n-1}|}{(2\pi)^n} \int_0^\infty \frac{\cos(tr)}{t} 
 \frac{\ud}{\ud r} \left[ r^{n-2} |\widehat{\sigma_2}(re_1)|^2 \right] \ud r \\[.4cm] \ds 
&=& - \ds\frac{|\mathbb{S}^{n-1}|}{(2\pi)^n} \int_0^\infty \frac{\sin(tr)}{t^2} 
 \frac{\ud^2}{\ud r^2} \left[ r^{n-2} |\widehat{\sigma_2}(re_1)|^2 \right] \ud r. 
\end{array} \]
Hence, we can estimate as follows 
\[ |q(t)| \leq \frac{K}{t^2} \quad \mbox{ with } \quad K= \frac{|\mathbb{S}^{n-1}|}{(2\pi)^n} 
 \int_0^\infty \left|  \frac{\ud^2}{\ud u^2} \left[ r^{n-2} |\widehat{\sigma_2}(re_1)|^2 \right] \right|\ud r<\infty  \]
which proves $q\in L^1([0,\infty))$.

Next, we compute the integral of $q$. For $M>0$ we get:
\[\begin{array}{lll}\ds \int_0^M q(t) \ud t 
&=& \ds\frac{1}{(2\pi)^n} \int_{\mathbb R^n} \left(\int_0^M  \ds \frac{\sin(t |\xi|) }{|\xi|} \ud t\right) |\widehat{\sigma_2}(\xi)|^2 \ud \xi 
\\
[.4cm]
&=&\ds\frac{1}{(2\pi)^n} \int_{\mathbb R^n} \frac{1-\cos(M|\xi|)}{|\xi|^2} |\widehat{\sigma_2}(\xi)|^2 \ud \xi \\[.4cm] \ds 
&=&\ds\kappa  - \frac{|\mathbb{S}^{n-1}|}{(2\pi)^n} \int_0^\infty \cos(Mr) r^{n-3} |\widehat{\sigma_2}(r e_1)|^2 \ud r \\ [.4cm]\ds 
&=&\ds\kappa  - \frac{|\mathbb{S}^{n-1}|}{M(2\pi)^n} \int_0^\infty \sin(Mr) \frac{\ud}{\ud r} \left[ r^{n-3} |\widehat{\sigma_2}(r e_1)|^2 \right]\ud r. 
\end{array}\]
We conclude by letting $M$ tend to $\infty$.

Note that $\kappa  $ is infinite for $n=2$ since $\frac{|\sigma_2(\xi)|^2}{|\xi|^2}\sim_{\xi\to 0} \|\sigma_2\|_{L^1(\mathbb R^2)}^2\frac1{|\xi|^2}$  does not belong to 
$L^1(B(0,a))$ for any $a>0$. 
\end{Proof}

We turn to the proof of Theorem \ref{cvth1}.
Of course we have
\[
\ds\sup_{\epsilon>0}\|f_\epsilon(t,\cdot)\|_{L^1(\mathbb R^d\times\mathbb R^d)}= \ds\sup_{\epsilon>0}\|f_{0,\epsilon}\|_{L^1(\mathbb R^d\times\mathbb R^d)}:=M_0<\infty,
\]
and the $L^p$ norms
\[
\|f_\epsilon(t,\cdot)\|_{L^p(\mathbb R^d\times\mathbb R^d)}= \|f_{0,\epsilon}\|_{L^p(\mathbb R^d\times\mathbb R^d)}\]
are also bounded, for any $1\leq p\leq\infty$ by virtue of  \eqref{H9}.
Furthermore, the energy conservation yields
\[\begin{array}{l}
\mathcal E_\epsilon(t)=
\ds\int_{\mathbb R^d\times\mathbb R^d}\Big(\ds\frac{v^2}{2} +V+\Phi_\epsilon \Big)f_\epsilon\ud v\ud x
\\[.4cm]
\qquad\qquad+\ds\frac\epsilon2\ds\int_{\mathbb R^n\times\mathbb R^d}|\partial_t \Psi_\epsilon|^2\ud y\ud x 
+\ds\frac1{2}\ds\int_{\mathbb R^n\times\mathbb R^d}|\nabla_y \Psi_\epsilon|^2\ud y\ud x \leq \bar{\mathcal E}_0.
\end{array}\]
Let us set 
\[
 \mathcal E_{0,\epsilon}^{\mathrm{vib}}=
 \ds\frac\epsilon2\ds\int_{\mathbb R^n\times\mathbb R^d}| \Psi_{1,\epsilon}|^2\ud y\ud x 
+\ds\frac1{2}\ds\int_{\mathbb R^n\times\mathbb R^d}|\nabla_y \Psi_{0,\epsilon}|^2\ud y\ud x.
\]
As a consequence of \eqref{H1}  and  \eqref{H8}, $ \mathcal E_{0,\epsilon}^{\mathrm{vib}}$ is bounded uniformly with respect to $\epsilon$.
Owing to the standard  energy conservation for the free linear wave equation, we observe that
$\| \nabla_y \widetilde \Psi_{\epsilon} \|_{L^\infty(0,\infty;L^2(\R^d \times \R^n))} \leq (2 \mathcal E_{0,\epsilon}^{\mathrm{vib}} )^{1/2}.  $ 
Then Sobolev's embedding (mind the condition $n\geq 3$) allows us to deduce the following key estimate on $\widetilde \Psi_{\epsilon}$:
\begin{equation}
\label{Psi0}
\| \widetilde \Psi_{\epsilon}\|_{L^\infty(\R_+;L^2(\R^d ; L^{2n/(n-2)} (\R^n)))} \leq C \big( \mathcal E_{0,\epsilon}^{\mathrm{vib}}  \big)^{1/2}
\leq C \big( \bar{\mathcal E}_{0}  \big)^{1/2}
\end{equation}
Applying H\"older 
 inequalities, we are thus led to:
 \begin{equation}
\label{Phi0}
|\Phi_{0,\epsilon}(t,x)|\leq C \|\sigma_2\|_{L^{2n/(n+2)}(\mathbb R^n)} \|\sigma_1\|_{L^2(\mathbb R^d)} \big( \bar{\mathcal E}_{0}  \big)^{1/2},
\end{equation}
and similarly
\begin{equation}
\label{Phi0g}
|\nabla_x\Phi_{0,\epsilon}(t,x)|\leq C \|\sigma_2\|_{L^{2n/(n+2)}(\mathbb R^n)} \| \nabla_x\sigma_1\|_{L^2(\mathbb R^d)} \big( \bar{\mathcal E}_{0}  \big)^{1/2}.
\end{equation}
Concerning the asymptotic behavior, we shall use the following claim. It is not a direct consequence of these estimates and it will be justified later on.
\begin{lemma}\label{phi0to0}
Let $\chi\in C^\infty_c([0,\infty)\times \mathbb R^d\times \mathbb R^d)$.
Then, we have
\[\ds\lim_{\epsilon\to 0}\ds\int_0^\infty\ds\int_{\mathbb R^d\times \mathbb R^d} f_\epsilon\nabla_x\Phi_{0,\epsilon}\chi(t,x,v)\ud v\ud x\ud t=0.\]
\end{lemma}

The cornerstone of the proof of Theorem \ref{cvth1} is the estimate of the self--consistent potential. By virtue of \eqref{resc_pot}, for any $1\leq p\leq\infty$ we get
\[\begin{array}{lll}  \left\|\ds\frac1\epsilon \mathcal{L}_\epsilon(f_\epsilon)(t,\cdot)\right\|_{L^p(\mathbb R^d)} &\leq&
  \| \Sigma \|_{L^p(\mathbb R^d)} \| \rho_\epsilon \|_{L^\infty([0,\infty),L^1(\mathbb R^d))}\ds \int_0^\infty |q(s)| \ud s 
\\[.4cm]
&\leq& \| \Sigma \|_{L^p(\mathbb R^d)} M_0 \|q\|_{L^1([0,+\infty))} ,\end{array} \]
as well as 
\[ \left\|\ds\frac1\epsilon \nabla_x \mathcal{L}_\epsilon(f_\epsilon)(t,\cdot)\right\|_{L^p(\mathbb R^d)}\leq \| \nabla_x\Sigma \|_{L^p(\mathbb R^d)} M_0 \|q\|_{L^1([0,+\infty))} .\]
Let $\chi\in C^\infty_c(\mathbb R^d\times \mathbb R^d)$.
We have
\[\left|\ds\int_{\mathbb R^d\times \mathbb R^d} f_\epsilon (t,x,v)\chi(x,v)\ud v\ud x\right| \leq M_0\|\chi\|_{L^\infty(\mathbb R^d\times \mathbb R^d)}
\]
and 
\[\begin{array}{l}
\left|\ds\frac{\ud}{\ud t}\ds\int_{\mathbb R^d\times \mathbb R^d} f_\epsilon (t,x,v)\chi(x,v)\ud v\ud x\right|
\qquad \leq  M_0
 \| v\cdot\nabla \chi - \nabla V \cdot \nabla_v \chi \|_{L^\infty(\mathbb R^d\times \mathbb R^d)}
\\[.4cm]
\qquad
+\Big( \|q\|_{L^1([0,+\infty))}  \|\nabla_x \Sigma \|_{L^\infty(\mathbb R^d)}M_0^2 + CM_0 \|\sigma_2\|_{L^{2n/(n+2)}(\mathbb R^n)} \| \nabla_x\sigma_1\|_{L^2(\mathbb R^d)} \big( \bar{\mathcal E}_{0}  \big)^{1/2} 
\Big)
\\[.4cm]
\hspace*{10cm}\times \|\nabla_v\chi\|_{L^\infty(\mathbb R^d\times \mathbb R^d)}.
\end{array} \]
Reproducing arguments detailed in the previous Section, we deduce that we can assume, possibly at the price of extracting a subsequence, that
\[\ds\lim_{\epsilon\rightarrow 0}\ds\int_{\mathbb R^d\times \mathbb R^d} f_\epsilon (t,x,v)\chi(x,v)\ud v\ud x
=\ds\int_{\mathbb R^d\times \mathbb R^d} f (t,x,v)\chi(x,v)\ud v\ud x
\] 
holds for any $\chi \in L^{p'}
(\mathbb R^d\times \mathbb R^d)$ uniformly on $[0,T]$, $0<T<\infty$, with $f\in C([0,T];L^p
(\mathbb R^d\times \mathbb R^d)-\text{weak})$, $1< p<\infty$, $1/p+1/p'=1$.

Next, we establish the tightness of $\big(f_\epsilon\big)_{\epsilon>0}$ with respect to the velocity variable, which will be necessary to show that the macroscopic density
$\rho_\epsilon$  passes to the limit.
Since  $\Phi_{0,\epsilon}$ and $\frac1\epsilon \mathcal{L}_\epsilon(f_\epsilon)$ are uniformly bounded 
and $V\geq 0$, we infer from  the energy conservation
the estimate
\[\begin{array}{l}\ds\int_{\mathbb R^d\times\mathbb R^d}\frac{|v|^2}{2} f_\epsilon(t,x,v)  \ud v \ud x
\\[.4cm]
\qquad\leq
 \bar{\mathcal{E}}_0 +  \|q\|_{L^1([0,+\infty))}  \| \Sigma \|_{L^\infty(\mathbb R^d)}M_0^2+CM_0 \|\sigma_2\|_{L^{2n/(n+2)}(\mathbb R^n)} \|\sigma_1\|_{L^2(\mathbb R^d)} \big( \bar{\mathcal E}_{0}  \big)^{1/2}
.\end{array}\]
Hence, we can check that $\rho_\epsilon(t,x)=\int_{\mathbb R^d}f_\epsilon(t,x,v)\ud v\ud x$
satisfies  \begin{equation}\label{cvrhoVl} \lim_{\epsilon\rightarrow0} \int_{\mathbb R^d} \rho_\epsilon(t,x) \chi(x) \ud x = \int_{\mathbb R^d} \rho(t,x) \chi(x) \ud x \end{equation}
for any $\chi\in C_0(\mathbb R^d)$, with $\rho(t,x)=\int_{\mathbb R^d} f(t,x,v)\ud v$.
As a matter of fact, we note that  \eqref{H1} and \eqref{cvrhoVl}
imply 
\begin{equation}\label{cvconv}
\ds\lim_{\epsilon\to 0} \nabla_x\Sigma\ast \rho_\epsilon(t,x)=\nabla_x\Sigma\ast\rho(t,x)\quad\text{for any $(t,x)\in [0,T]\times\mathbb R^d$.}
\end{equation}
Furthermore, we have
\[  |D^2_x (\Sigma\ast\rho_\epsilon)(t,x) | \leq M_0 \ \| \Sigma \|_{W^{2,\infty}(\mathbb R^d)}, \]
and, by using mass conservation and the Cauchy-Schwarz inquality,
\[\begin{array}{l}
|\partial_t ( \nabla_x\Sigma\ast \rho_\epsilon)(t,x)|=\left|\ds\int_{\mathbb R^d}D^2_x\Sigma(x-y)\left(\ds\int_{\mathbb R^d} vf_\epsilon(t,y,v)\ud v\right)\ud y\right|
\\[.4cm]\quad\leq 
\| \Sigma \|_{W^{2,\infty}(\mathbb R^d)}\left(\ds\int_{\mathbb R^d\times\mathbb R^d} f_\epsilon\ud v\ud x\right)^{1/2}\left(\ds\int_{\mathbb R^d\times\mathbb R^d}v^2 f_\epsilon\ud v\ud x\right)^{1/2}
\\[.4cm]\quad\leq  \| \Sigma \|_{W^{2,\infty}(\mathbb R^d)}
\ \sqrt{2M_0}\ \Big(\bar{\mathcal{E}}_0 +  \|q\|_{L^1([0,+\infty))}  \| \Sigma \|_{L^\infty(\mathbb R^d)}M_0^2 
\\[.4cm]\qquad\qquad\qquad\qquad+CM_0 \|\sigma_2\|_{L^{2n/(n+2)}(\mathbb R^n)} \|\sigma_1\|_{L^2(\mathbb R^d)} \big( \bar{\mathcal E}_{0}  \big)^{1/2}\Big)^{1/2}
.\end{array}\]
Therefore convergence \eqref{cvconv}
holds uniformly on any compact set  of $[0,\infty)\times\mathbb R^d$.
%

We turn to examine the convergence of  $\frac1\epsilon\nabla_x \mathcal{L}_\epsilon (f_\epsilon)$ to  $\kappa \nabla_x\Sigma \ast \rho$.
We have
\[\begin{array}{l}
\Big|\ds\frac1\epsilon\nabla_x \mathcal{L}_\epsilon (f_\epsilon)(t,x)- \kappa \nabla_x\Sigma \ast \rho(t,x)\Big|
\\[.4cm]
\qquad= 
 \left| \ds\int_0^{t/\sqrt\epsilon} \nabla_x\Sigma \ast \rho_\epsilon(t-s\sqrt\epsilon,x) q(s) \ud s - \kappa  \nabla_x\Sigma \ast \rho(t,x) \right| \\[.4cm]
\qquad \leq\ds 
 \left| \ds\int_0^{t/\sqrt\epsilon} \left(\nabla_x\Sigma \ast \rho_\epsilon(t-s\sqrt\epsilon,x)-\nabla_x\Sigma \ast \rho(t,x) \right) q(s) \ud s \right| 
 \\
 [.4cm]\qquad\qquad
+ 
 \left|\ds \int_{t/\sqrt\epsilon} ^\infty q(s) \ud s \right| \| \nabla_x\Sigma  \ast \rho \|_{L^\infty((0,\infty)\times\mathbb R^d)}
 \\[.4cm]
 \qquad\leq
\ds \int_0^{t/\sqrt\epsilon} |(\nabla_x\Sigma \ast \rho_\epsilon- \nabla_x\Sigma \ast \rho)(t-s\sqrt\epsilon,x)| \  |q(s)| \ud s 
\\[.4cm]\qquad\qquad
+ \ds 
 \int_0^{t/\sqrt\epsilon}  |\nabla_x\Sigma \ast \rho(t-s\sqrt\epsilon,x) - \nabla_x\Sigma \ast \rho(t,x)|\  |q(s)| \ud s
 \\[.4cm]\qquad\qquad
+ 
\ds \int_{t/\sqrt\epsilon} ^\infty |q(s)| \ud s \  \| \nabla_x\Sigma  \ast \rho \|_{L^\infty((0,\infty)\times\mathbb R^d)}
 .
\end{array} \]
Let us denote by $\mathrm I_\epsilon(t,x)$, $\mathrm {II}_\epsilon(t,x)$, $\mathrm {III}_\epsilon(t)$, the three terms of the right hand side.
Firstly, for any $t>0$, $\mathrm {III}_\epsilon(t)$ tends to 0 as $\epsilon\to 0$, and it is dominated by $\kappa  \|\Sigma\|_{W^{1,\infty}(\mathbb R^d)}M_0$.
Secondly, for any $0<T<\infty$ and any compact set $K\subset\mathbb R^d$, when $(t,x)$ lies in $[0,T]\times K$, we
 can estimate $$ | \mathrm {I}_\epsilon(t,x) | \leq \|\nabla_x\Sigma \ast \rho_\epsilon- \nabla_x\Sigma  \ast \rho\|_{L^\infty([0,T]\times K)}
 \|q\|_{L^1([0,\infty)}$$
 which also goes to 0 as $\epsilon\to 0$.
 Eventually, still considering $(t,x)\in [0,T]\times K$,
 we write
 \[\begin{array}{lll}
 |\mathrm{II}_\epsilon(t,x)| 
  \leq\ds 
 \int_0^{t/\sqrt\epsilon}  \ds\sup_{z\in K}| \nabla_x\Sigma \ast \rho(t-s\sqrt\epsilon,z)-\nabla_x\Sigma \ast \rho(t,z)|\ |q(s)|\ud s.
\end{array} 
 \]
 By using the Lebesgue theorem, we justify that it tends to 0
  as $\epsilon\to 0$ since $(t,x)\mapsto \nabla_x\Sigma \ast \rho(t,x)$ is uniformly continuous over any compact set, the integrand is dominated by 
  $2\|\Sigma\|_{W^{1,\infty}(\mathbb R^d)} M_0|q(s)|$, and $q\in L^1([0,\infty))$.
Therefore, for any $0<t<T<\infty$ and any compact set $K\subset\mathbb R^d$, 
$$
\ds\sup_{x\in K}\Big|\frac1\epsilon\nabla_x \mathcal{L}_\epsilon (f_\epsilon)-\kappa \nabla_x\Sigma \ast \rho\Big|(t,x)
\xrightarrow[\epsilon\to 0] {}0,
$$
and this quantity is bounded uniformly with respect to $0\leq t\leq T<\infty$ and $\epsilon>0$.

We 
go back to the weak formulation of  \eqref{limite1}. Let $\chi\in C^\infty_c([0,\infty)\times \mathbb R^d \times \mathbb R^d)$.
We suppose that $\mathrm {supp}(\chi)\subset [0,T]\times \bar B(0,M)\times \bar B(0,M)$.
We have
\[ \begin{array}{l}
-\ds \int_{\mathbb R^d \times \mathbb R^d} f_{0,\epsilon} \chi(0,x,v) \ud v \ud x-
 \int_0^\infty\int _ {\mathbb R^d \times \mathbb R^d}f_\epsilon \partial_t  \chi \ud v \ud x\ud t 
 \\
 [.4cm]\qquad\ds
 -
  \int_0^\infty\int _ {\mathbb R^d \times \mathbb R^d} f_\epsilon v\cdot\nabla_x \chi \ud v \ud x\ud t 
+  \int_0^\infty\int _ {\mathbb R^d \times \mathbb R^d}\int f_\epsilon \nabla_v \chi \cdot \nabla_x   (V+\Phi_{0,\epsilon})
\ud v \ud x\ud t
\\
 [.4cm]\qquad\qquad\ds= \int_0^\infty\int _ {\mathbb R^d \times \mathbb R^d}
f_\epsilon\ \nabla_x \ds\frac1\epsilon \mathcal{L}_\epsilon(f_\epsilon)\cdot\nabla_v \chi \ud v \ud x\ud t .\end{array}\]
Obviously, there is no  difficulty with the linear terms of the left hand side.
For the non linear term we proceed as follows:
\[\begin{array}{l}
\ds\int_0^\infty\int _ {\mathbb R^d \times \mathbb R^d}
f_\epsilon\ \nabla_x \ds\frac1\epsilon \mathcal{L}_\epsilon(f_\epsilon)\cdot\nabla_v \chi \ud v \ud x\ud t
-
\int_0^\infty\int _ {\mathbb R^d \times \mathbb R^d}
f\ \kappa  \nabla_x\Sigma\ast\rho \cdot\nabla_v \chi \ud v \ud x\ud t
\\[.4cm]
\qquad\ds
=\int_0^\infty\int _ {\mathbb R^d \times \mathbb R^d}
f_\epsilon\ \Big(\nabla_x \ds\frac1\epsilon \mathcal{L}_\epsilon(f_\epsilon)-\kappa  \nabla_x\Sigma\ast\rho \Big)\cdot\nabla_v \chi \ud v \ud x\ud t
\\[.4cm]
\qquad\qquad+
\ds\int_0^\infty\int _ {\mathbb R^d \times \mathbb R^d}
(f_\epsilon-f)\ \kappa  \nabla_x\Sigma\ast\rho \cdot\nabla_v \chi \ud v \ud x\ud t.
\end{array}\]
 The last term directly passes to the limit.
 The first integral in the right hand side  is dominated by
 \[
 M_0\|\nabla_v\chi\|_{L^\infty([0,\infty)\times\mathbb R^d\times\mathbb R^d)}\ 
  \ds\int_0^T 
 \ds\sup_{y\in \bar B(0,M)} \Big|
 \nabla_x \ds\frac1\epsilon \mathcal{L}_\epsilon(f_\epsilon)-\kappa  \nabla_x\Sigma\ast\rho\Big|(t,y)\ud t.
 \]
 We conclude by a mere application of the Lebesgue Theorem.

If the initial data $f_{0,\epsilon}$ converge strongly to $f_0$ in $L^p(\mathbb R^d\times\mathbb R^d)$, 
the nature of the convergence of $f_\epsilon$ to $f$ can be  improved
by applying general stability results for transport equations, see \cite[Th. II.4 \& Th. II.5]{DPL}, or \cite[Th. VI.1.9]{BF}.
\\

\noindent
{\bf Proof of Lemma \ref{phi0to0}}
As a matter of fact, the variable $x\in\mathbb R^d$ 
just appears as a parameter for the wave equation, and $\Upsilon_\epsilon(t,x,y)=(\sigma_1\ast\widetilde\Psi_\epsilon(t,\cdot,y))(x)$
solves the linear wave equation
\[\epsilon\partial^2_{tt}\Upsilon_\epsilon-\Delta_y\Upsilon_\epsilon=0,\]
with the data
\[\Upsilon_\epsilon(0,x,y)=\sigma_1\ast\Psi_{0,\epsilon}(x,y),\qquad\partial_t\Upsilon_\epsilon(0,x,y)=\sigma_1\ast\Psi_{1,\epsilon}(x,y).\]
The parameter  $x$ being fixed, we appeal to the Strichartz estimate, see  \cite[Corollary 1.3]{MT} or \cite[Theorem 4.2, for the case $n=3$]{Sog},
\[
\ds\frac{1}{\epsilon^{1/(2p)}}\left(\ds\int_0^\infty\left(\ds\int_{\mathbb R^n}|\Upsilon_\epsilon(t,x,y)|^q\ud y\right)^{p/q}\ud t\right)^{1/p}\leq  C \sqrt{\mathscr E^{\mathrm{vib}}_{1,\epsilon}(x)}\]
where 
we set \[
\mathscr E^{\mathrm{vib}}_{1,\epsilon}(x)=\epsilon\ds\int_{\mathbb R^n}|\sigma_1\ast\Psi_{1,\epsilon}(x,y)|^2\ud y+ \ds\int_{\mathbb R^n}|\sigma_1\ast\nabla_y\Psi_{0,\epsilon}(x,y)|^2\ud y
.\]
(That $\frac{1}{\epsilon^{1/(2p)}}$ appears in the inequality can be checked by changing variables and observing that $\Upsilon_\epsilon(t\sqrt \epsilon ,x,y)$ satisfies
the wave equation with speed equals to 1 and data $(\sigma_1\ast\Psi_{0\epsilon},\sqrt\epsilon\sigma_1\ast\Psi_{1,\epsilon})$.)
This inequality holds for admissible exponents:
\[2\leq p\leq q\leq\infty,\quad
\ds\frac1p+\ds\frac n q=\ds\frac n2-1,\quad
\ds\frac2p+\ds\frac{n-1}q\leq \ds\frac {n-1}2,\quad
(p,q,n)\neq (2,\infty,3).\]
Observe that 
\[
\ds\int_{\mathbb R^d} \mathscr E^{\mathrm{vib}}_{1,\epsilon}(x) \ud x
\leq \| \sigma_1\|_{L^1(\mathbb R^d)}\ \mathscr E^{\mathrm{vib}}_{0,\epsilon}
\leq  \| \sigma_1\|_{L^1(\mathbb R^d)}\ \bar{\mathscr E}_{0}.
\]
%
%
It follows that
\[\ds\int_{\mathbb R^d}
\left(\ds\int_0^\infty\left(\ds\int_{\mathbb R^n}|\Upsilon_\epsilon(t,x,y)|^q\ud y\right)^{p/q}\ud t
\right)^{2/p}\ud x
\leq 
 C^2\| \sigma_1\|_{L^1(\mathbb R^d)}\ \bar{\mathscr E}_{0}\epsilon^{1/p}
\xrightarrow[\epsilon\to 0]{}0.\]
A similar reasoning applies to $\nabla_x\Upsilon_\epsilon$ with $\nabla_x\sigma_1$ replacing $\sigma_1$.
Let $\chi\in C^\infty_c([0,\infty)\times\mathbb R^d\times\mathbb R^d)$. We suppose that $\mathrm{supp}(\chi)\subset \{0\leq t \leq M,\ |x|\leq M,\ |v|\leq M\}$ for some $0<M<\infty$.
We are left with the task of estimating
\[
\ds\int_0^\infty\ds\int_{\mathbb R^d\times\mathbb R^d} f_\epsilon \nabla_x \Phi_{0,\epsilon}\chi(t,x,v)\ud v\ud x\ud t
=
\ds\int_0^\infty\ds\int_{\mathbb R^d} R_\epsilon(t,x) \nabla_x \Phi_{0,\epsilon}(t,x)\ud x\ud t
\]
where we have set 
\[R_\epsilon(t,x)=\ds\int_{\mathbb R^d}f_\epsilon\chi(t,x,v)\ud v.\]
With the standard notation $1/p+1/p'=1$, using H\"older's inequality twice, we get
\[\begin{array}{l}
\left|\ds\int_0^\infty\ds\int_{\mathbb R^d\times\mathbb R^d} f_\epsilon \nabla_x \Phi_{0,\epsilon}\chi(t,x,v)\ud v\ud x\ud t
\right|
\\[.4cm]
\quad\leq 
\left(\ds\int_{\mathbb R^d}
\left(
\ds\int_0^\infty |R_\epsilon(t,x)|^{p'}\ud t\right)^{2/p'}
\ud x\right)^{1/2}
\left(\ds\int_{\mathbb R^d}
\left(
\ds\int_0^\infty |\nabla_x\Phi_{0,\epsilon}(t,x)|^{p}\ud t\right)^{2/p}
\ud x\right)^{1/2}.
\end{array}\]
We readily obtain 
\[\begin{array}{lll}
\left(\ds\int_{\mathbb R^d}
\left(
\ds\int_0^\infty |R_\epsilon(t,x)|^{p'}\ud t\right)^{2/p'}
\ud x\right)^{1/2}
&\leq& M^{d+d/2+1/p'}\|f_\epsilon \chi \|_{L^\infty((0,\infty)\times \mathbb R^d\times\mathbb R^d)} \\
&\leq &M^{d+d/2+1/p'}\|f_{0,\epsilon}\|_{L^\infty(\mathbb R^d\times\mathbb R^d)} 
\| \chi \|_{L^\infty((0,\infty)\times \mathbb R^d\times\mathbb R^d)}
\end{array}\]
which is thus 
bounded uniformly with respect to $\epsilon>0$.
Furthermore, with $1/q+1/q'=1$, we have
\[\begin{array}{l}
\ds\int_{\mathbb R^d}\left(\ds\int_0^\infty |\nabla_x\Phi_{0,\epsilon}(t,x)|^{p}\ud t\right)^{2/p}
\ud x
=
\ds\int_{\mathbb R^d}\left(\ds\int_0^\infty \Big|
\ds\int_{\mathbb R^n}\sigma_2(y)
\nabla_x\Upsilon_{\epsilon}(t,x,y)\ud y
\Big|^{p}\ud t\right)^{2/p}\ud x
\\
\qquad
\leq
\|\sigma_2\|_{L^{q'}(\mathbb R^d)}
\ds\int_{\mathbb R^d}\left(\ds\int_0^\infty \Big|
\ds\int_{\mathbb R^n}
|\nabla_x\Upsilon_{\epsilon}(t,x,y)|^q\ud y
\Big|^{p/q}\ud t\right)^{2/p}\ud x
\end{array}
\]
which tends to 0 like $\epsilon^{1/p}$.
\QED

\subsection{Convergence to the Vlasov--Poisson system}

The existence theory for the Vlasov--Poisson system dates back to  \cite{Ars}; an overview of the features of both the   repulsive or attractive cases can be found in  the lecture notes \cite{Bouc}.
The following statements are classical tools of this analysis, that will be useful for our purposes as well.
 \begin{lemma}[Interpolation estimates]\label{interp}
Let $f \in L^1 \cap L^\infty(\mathbb R^d \times \mathbb R^d)$ be such that  $ |v|^m f\in L^1(\mathbb R^d \times \mathbb R^d)$.
 Then $\rho=\int_{\mathbb R^d} f\ud v$ lies in  $L^{(m+d)/d}(\mathbb R^d)$ with
\[ \| \rho \|_{L^{(d+m)/d}(\mathbb R^d)} \leq C(m,d) \|f\|_{L^\infty}^{m/(d+m)(\mathbb R^d)} \left( \int  |v|^m f\ud v \ud x\right)^{d/(d+m)}.  \]
where $C(m,d) =2|B(0,1)|^{m/(m+d)}$.
\end{lemma}

\begin{lemma}[Hardy-Littlewood-Sobolev inequality]\label{HLS}
Let $1<p,r<\infty$ and $0<\lambda<d$.
Assume $1/p+1/r=2-\lambda/d$.
There exists a constant $C>0$ such that  
for any $f\in L^p(\mathbb R^d)$ and $g\in L^r(\mathbb R^d)$
 we have
\[
\Big|\ds\int_{\mathbb R^d\times\mathbb R^d} \ds\frac{f(x)g(y)}{|x-y|^\lambda}\ud y\ud x\Big|\leq C\|f\|_{L^p(\mathbb R^d)}\|g\|_{L^r(\mathbb R^d)}.
\]
\end{lemma}

We refer the reader to 
\cite[Lemma 3.4]{Bouc} and 
\cite[Th. 4.3]{LL}, respectively, for further details.
Next, we check the convergence of the approximate kernel defined by $\sigma_{1,\epsilon}$.

\begin{lemma}\label{kernel}
Let $d\geq 3$. For any  
$ d/(d-1)<q<\infty$, we have:
\[ \left\|  \nabla \left( \frac{C_d\theta_\epsilon}{|\cdot|^{d-1}} \ast\frac{C_d\theta_\epsilon}{|\cdot|^{d-1}} \right) (x)+ (d-2)\frac{x}{|\mathbb{S}^{d-1}| |x|^d} \right\| _{L^q(\mathbb R^d)} 
\xrightarrow[\epsilon\to0]{} 0. \]
\end{lemma}

\begin{Proof}
We remind the reader that the convolution by $|x|^{1-d}$ is associated to the Fourier transform of the operator with symbol $1/|\xi|$, see \cite[Th. 5.9]{LL}.
The convolution of radially symmetric functions is radially symmetric too. For  $d \geq 3$, we compute as follows
\[\begin{array}{lll}\left( \ds\frac{1}{|\cdot|^{d-1}} \ast\frac{1}{|\cdot|^{d-1}} \right)(x) 
&=& \ds\int_{\mathbb R^d} \frac{\ud y}{|y|^{d-1}|x-y|^{d-1}}
\\[.4cm]
&=&\ds \int_{\mathbb R^d} \frac{|x|^d \ud y}{|x|^{d-1}|e_1 -y|^{d-1} |x|^{d-1} |y|^{d-1}}
= \frac{1}{|\mathbb{S}^{d-1}| \ C_d^2\ |x|^{d-2}}. \end{array}\]
Differentiating yields
\[ \nabla \left( \frac{C_d}{|\cdot|^{d-1}} \ast\frac{C_d}{|\cdot|^{d-1}} \right)(x) = -\frac{d-2}{ |\mathbb{S}|^{d-1}}\  \frac{x}{ |x|^d}. \]
Hence, we can write
\[ \begin{array}{lll}\mathscr O_\epsilon(x)&:=&\ds\nabla \left( \frac{C_d\theta_\epsilon}{|\cdot|^{d-1}} \ast\frac{C_d\theta_\epsilon}{|\cdot|^{d-1}} \right)(x) +\frac{(d-2)x}{|\mathbb{S}^{d-1}| |x|^d}  \ds 
\\[.4cm]
&=& \ds C_d^2 \nabla \left( \frac{\theta_\epsilon +1}{|\cdot|^{d-1}} \ast  \frac{\theta_\epsilon -1}{|\cdot|^{d-1}} \right)(x) \\[.4cm] \ds 
&=&  \ds C_d^2 \frac{\theta_\epsilon+1}{|\cdot|^{d-1}} \ast \left( \frac{\nabla \theta_\epsilon }{|\cdot|^{d-1}}  + (1-d) \frac{(\theta_\epsilon -1)\cdot}{|\cdot|^{d+1}}  \right)(x).
\end{array} \]
Let $p>1$.
On the one hand, we
have
\[ \begin{array}{lll}
\ds\left\| \frac{\nabla \theta_\epsilon}{|\cdot|^{d-1}}  \right\|^p_{L^p(\mathbb R^d)} 
&=&\ds \int_{\mathbb R^d} \frac{|\nabla \theta_\epsilon(x)|^p}{|x|^{p(d-1)}} \ud x
\\[.4cm]
&\leq&(\sqrt\epsilon)^p\ \|\nabla \theta\|_{L^\infty(\mathbb R^d)}^p\ \ds \int_{1 \leq \sqrt\epsilon |x| \leq 2}  \frac{\ud x}{|x|^{p(d-1)}}
\\[.4cm]
&\leq& (\sqrt\epsilon)^{d(p-1)}
 \|\nabla \theta\|_{L^\infty(\mathbb R^d)}^p\ds \int_{1 \leq |x| \leq 2}  \frac{\ud x}{|x|^{p(d-1)}}
. \end{array} \]
On the other hand, we get
\[\ds\int_{\mathbb R^d} \Big|\frac{(\theta_\epsilon(x)-1)x }{|x|^{d+1}}\Big|^p \ud x 
\leq \int_{\sqrt\epsilon |x|\geq 1} \frac{\ud x}{|x|^{pd}}
= (\sqrt\epsilon)^{d(p-1)}\left( \int_{|x|\geq 1}   \frac{\ud x}{|x|^{pd}} \right).  \]
Accordingly, the following estimate holds:
\begin{equation}
\label{sigma1}
\left\|  \frac{\nabla \theta_\epsilon }{|\cdot|^{d-1}}  + (1-d) \frac{(\theta_\epsilon -1)\cdot}{|\cdot|^{d+1}} \right\|_{L^p} 
\leq C \epsilon^{d(p-1)/(2p)},
\end{equation}         
where $C>0$ depends on $p$ and $d$ only.
Finally we remark that $ 0 \leq \frac{\theta_\epsilon(x) +1}{|x|^{d-1}} \leq \frac{2}{|x|^{d-1}}$.
By coming back to Lemma \ref{HLS}, we deduce that
there exists a constant $\tilde C>0$ such that
\[
\Big|\ds\int_{\mathbb R^d}\mathscr O_\epsilon(x)g(x)\ud x\Big|\leq \tilde C\|g\|_{L^r(\mathbb R^d)}\ (\sqrt\epsilon)^{d(p-1)/p}
\]
holds for any $g\in L^r(\mathbb R^d)$, with $ 1/r=(d+1)/d-1/p>1/d$, $r>1$.
Therefore, by duality, it means that $\mathscr O_\epsilon$ converges to 0 in $L^q(\mathbb R^d)$ for any $ d/(d-1)<q<\infty$.
\end{Proof}

\begin{ProofOf}{Theorem \ref{cvth2}}
From now on, we restrict to the case of space dimension $d=3$.
Compared to the previous Section,  additional difficulties come from the dependence of the form function $\sigma_1$ with respect to $\epsilon$
so that  deducing uniform estimates from the energy conservation is not direct.

\noindent{\it Step 1. Establishing uniform estimates.}

We start by observing that $f_\epsilon$ is bounded in $L^\infty(0,\infty;L^p(\mathbb R^3\times\mathbb R^3))$ for any $1\leq p\leq \infty$, since
\[
\|f_\epsilon(t,\cdot)\|_{L^p(\mathbb R^3\times\mathbb R^3)}= \|f_{0,\epsilon}\|_{L^p(\mathbb R^3\times\mathbb R^3)}.\]
Next, the energy conservation becomes
 \[\begin{array}{lll}
 \mathcal{E}_\epsilon(t)&=& \ds\frac{\epsilon}{2} \int_{\mathbb R^3\times \mathbb R^n} \left| \partial_t \Psi_\epsilon(t,x,y) \right|^2 \ud y \ud x 
+ \frac{1}{2} \int_{\mathbb R^3\times \mathbb R^n} \left| \nabla_y \Psi_\epsilon(t,x,y) \right|^2 \ud y \ud x 
\\
&&+\ds \int_{\mathbb R^3\times \mathbb R^3} f_\epsilon(t,x,v) \left( \frac{|v|^2}{2} + V(x) + \Phi_\epsilon(t,x) \right) \ud v \ud x 
\\
&=&\mathcal{E}_\epsilon(0)\leq \bar{\mathcal E}_0.\end{array}\]
Let us study the coupling term:
\[  \ds\int_{\mathbb R^3\times \mathbb R^3} f_{\epsilon}(t,x,v) \Phi_\epsilon(t,x) \ud v \ud x \ds 
=\ds \int_{\mathbb R^3} \rho_{\epsilon}(t,x) \Phi_\epsilon(t,x) \ud x=\mathrm{S}_\epsilon(t)+\mathrm{T}_\epsilon(t)
\]
where we have set
\[\begin{array}{lll}
\mathrm{S}_\epsilon(t)
&=&-\ds\frac{1}{\epsilon}\ds \int_{\mathbb R^3}\rho_{\epsilon}\mathcal{L}_\epsilon(f_{\epsilon})(t,x)\ud x
\\
&=&
- \ds \int_{\mathbb R^3}\Big(\sigma_{1,\epsilon} \ast \sigma_{1,\epsilon} \ast\ds \int_0^{t/\sqrt\epsilon} q(s) \rho_{\epsilon} (t-s\sqrt\epsilon,\cdot) \ud s\Big)(x)\rho_\epsilon(t,x)\ud x
\\
&=&- \ds \int_{\mathbb R^3}\Big( \sigma_{1,\epsilon} \ast\ds \int_0^{t/\sqrt\epsilon} q(s) \rho_{\epsilon} (t-s\sqrt\epsilon,\cdot) \ud s\Big)(x)\ \sigma_{1,\epsilon} \ast\rho_\epsilon(t,x)\ud x
\end{array}\]
and 
\[
\mathrm{T}_\epsilon(t)=\ds\int_{\mathbb R^3} \rho_{\epsilon} \Phi_{0,\epsilon}(t,x)\ud x,
\qquad \Phi_{0,\epsilon}(t,x) = \left(\sigma_{1,\epsilon} \ast \int_{\R^n} \widetilde \Psi_{\epsilon}(t,\cdot,y) \sigma_2(y) \ud y\right)(x). \]
Like in the previous Section, $\widetilde \Psi_{\epsilon}$ stands for the  solution of the free linear wave equation with wave speed $1/\epsilon$ and  initial data $\Psi_{0,\epsilon}$ and 
$\Psi_{1,\epsilon}$.
Firstly, we 
 establish a bound for
\[
|\mathrm {S}_\epsilon(t)|\leq 
\|q\|_{L^1([0,\infty))}\|\sigma_{1,\epsilon} \ast \rho_{\epsilon}\|_{L^\infty(0,t;L^2(\mathbb R^3))}^2.
\]
However, Lemma \ref{HLS} yields
\[
\|\sigma_{1,\epsilon} \ast \rho_{\epsilon}\|_{L^2(\mathbb R^3)}
=C_d^2 \left\|\ds \frac{\theta_\epsilon}{|\cdot|^{2}} \ast \delta_\epsilon \ast \rho_{\epsilon}\right\|_{L^2(\mathbb R^3)}
\leq C\| \rho_\epsilon\|_{L^{6/5}(\mathbb R^3)}.
\]
Let us set  
\[\mathcal E^{\mathrm{kin}}_{\epsilon}(t)=\ds\int_{\mathbb R^3\times\mathbb R^3}|v|^2 f_\epsilon(t,x,v)\ud v\ud x\]
for the particle kinetic energy.
Lemma \ref{interp} leads to
\begin{equation}
\label{rho}
\| \rho_{\epsilon} \|_{L^{5/3}(\mathbb R^3)}  \leq C(2,3) \|f_{\epsilon} \|_{L^\infty(\mathbb R^3\times\mathbb R^3)}^{2/5}
\big(\mathcal E^{\mathrm{kin}}_{\epsilon}\big)^{3/5} 
\end{equation}
The H\"older inequality allows us to estimate 
$ \| \rho_{\epsilon}\|_{L^{6/5}(\mathbb R^3)}  \leq \| \rho_{\epsilon} \|_{L^1(\mathbb R^3)}^{7/12}  \| \rho_{\epsilon} \|_{L^{5/3}(\mathbb R^3)}^{5/12} $.
Combining these inequalities, we arrive at 
\begin{equation}\label{conv} \| \sigma_{1,\epsilon} \ast \rho_{\epsilon} \|_{L^2(\mathbb R^3)} 
\leq C \big(\mathcal E^{\mathrm{kin}}_{\epsilon}\big)^{1/4},\end{equation}
for a certain constant $C>0$, which does not depend on $\epsilon$.
Therefore, we obtain 
\[ |\mathrm {S}_\epsilon(t)| \leq C^2 \|q\|_{L^1([0,\infty))} \| \mathcal E^{\mathrm{kin}}_{\epsilon} \|_{L^\infty([0,t])}^{1/2}.\]
Secondly, we  estimate the  term involving $\Phi_{0,\epsilon}$:
\[ \mathrm {T}_\epsilon(t)
= \int_{\R^d\times \R^N} (\rho_\epsilon \ast \sigma_{1,\epsilon} )(t,x) \widetilde \Psi_\epsilon(t,x,y) \sigma_2(y) \ud y \]
 is dominated by
\[ \|\sigma_{1,\epsilon} \ast \rho_{\epsilon}\|_{L^\infty(0,t;L^2(\mathbb R^3))} 
\| \widetilde \Psi_\epsilon \|_{L^\infty(\R_+ ; L^2(\R^d ; L^{2n/(n-2)}(\R^n)))} 
\| \sigma_2 \|_{L^{2n/(n+2)}(\R^n)}. \]
Using \eqref{Psi0} and \eqref{conv}, we get
\[| \mathrm {T}_\epsilon(t)| \leq C' \big(\mathcal E^{\mathrm{kin}}_{\epsilon}(t)\big)^{1/4} \big( \mathcal E_{0,\epsilon}^{\mathrm{vib}} \big)^{1/2} \]
where the constant $C'>0$ does not  depend on $\epsilon$.
It remains to discuss how \eqref{H7}--\eqref{H8} implies a uniform estimate on the initial state.
Note that $\mathrm{S}_\epsilon(0)=0$.
Hence, by using \eqref{H8}, we are led to 
\[
\mathcal E_{0,\epsilon}^{\mathrm{vib}} + \frac12 \mathcal E^{\mathrm{kin}}_{\epsilon}(0)
\leq \mathcal E_\epsilon (0)+ |\mathrm {T}_\epsilon(0) |\leq 
 \bar{\mathcal{E}}_0 + C' \big(\mathcal E^{\mathrm{kin}}_{\epsilon}(0)\big)^{1/4} \big( \mathcal E_{0,\epsilon}^{\mathrm{vib}} \big)^{1/2}.
\]
It allows us to infer 
\[ \ds\sup_{0<\epsilon<1} \mathcal E^{\mathrm{kin}}_{\epsilon}(0)=\bar{\mathcal E}^{\mathrm{kin}}_0<\infty
,
\qquad 
\ds\sup_{0<\epsilon<1}  \mathcal E_{0,\epsilon}^{\mathrm{vib}} =\bar{\mathcal E}_{0}^{\mathrm{vib}}<\infty.
\]
Coming back to the energy conservation, with\eqref{H7}--\eqref{H8}  together with the estimates on $\mathrm T_\epsilon$ and $\mathrm S_\epsilon$, we deduce that
\[\ds
\frac12 \mathcal E^{\mathrm{kin}}_{\epsilon}(t)\leq
\bar{\mathcal{E}}_0 + C^2 \|q\|_{L^1([0,\infty))} \| \mathcal E^{\mathrm{kin}}_{\epsilon} \|_{L^\infty([0,t])}^{1/2} 
+
C' \big(\mathcal E^{\mathrm{kin}}_{\epsilon}(t)\big)^{1/4} \big( \bar{ \mathcal E}_{0,\epsilon}^{\mathrm{vib}} \big)^{1/2}  ,\]
holds, which, in turn, establishes the bound 
\[ \ds\sup_{0<\epsilon<1,\ t\geq 0} \mathcal E^{\mathrm{kin}}_{\epsilon}(t) =\bar{\mathcal E}^{\mathrm{kin}}<\infty. \] 
Going back to the interpolation inequalities, it follows that $\rho_\epsilon$ is bounded in $L^\infty(0,\infty;L^1\cap L^{5/3}(\mathbb R^3))$.
\\

\noindent
{\it Step 2. Passing to the limit.}

The kinetic equation can be rewritten
\[\partial_t f_{\epsilon}+v\cdot\nabla_x f_\epsilon -\nabla_x\Big(V +   \Phi_{0,\epsilon} -\ds\frac1\epsilon \mathcal L_\epsilon(f_\epsilon)\Big)\cdot\nabla_v f_\epsilon=0.
\]
We start by establishing that $ \nabla_v f_\epsilon \cdot \nabla_x \Phi_{0,\epsilon} =\nabla_v\cdot(
f_\epsilon  \nabla_x \Phi_{0,\epsilon})$ converges to $0$ at least in the sense of distributions.

 \begin{lemma}\label{phi0to0vp}
Let $\chi\in C^\infty_c([0,\infty)\times \mathbb R^d\times \mathbb R^d)$.
Then, we have
\[\ds\lim_{\epsilon\to 0}\ds\int_0^\infty\ds\int_{\mathbb R^d\times \mathbb R^d} f_\epsilon\nabla_x\Phi_{0,\epsilon} \chi(t,x,v)\ud v\ud x\ud t=0.\]
\end{lemma}

\begin{Proof}
It is convenient to split
\[\begin{array}{lll}
\Phi_{0,\epsilon}(t,x)&=&
\ds\int_{\mathbb R^n} \sigma_2(y) C_3\ds\frac{\theta_\epsilon}{|\cdot|^2}\ast\delta_\epsilon\ast \widetilde \Psi_\epsilon(t,x,y) \ud y
\\[.4cm]
&=&
\Phi_{0,\epsilon}^{\mathrm{main}}(t,x)+\Phi_{0,\epsilon}^{\mathrm{rem}}(t,x)
\end{array}\]
with
\[\begin{array}{l}
\Phi_{0,\epsilon}^{\mathrm{main}}(t,x)=
\ds\int_{\mathbb R^n} \sigma_2(y) C_3\ds\frac{1}{|\cdot|^2}\ast\delta_\epsilon\ast \widetilde \Psi_\epsilon(t,x,y) \ud y,
\\[.4cm]
\Phi_{0,\epsilon}^{\mathrm{rem}}(t,x)=
\ds\int_{\mathbb R^n} \sigma_2(y) C_3\ds\frac{\theta_\epsilon-1}{|\cdot|^2}\ast\delta_\epsilon\ast \widetilde \Psi_\epsilon(t,x,y) \ud y,
\end{array}
\]
and we remind the reader that $\widetilde \Psi_\epsilon(t,x,y) $ is the solution 
of the free wave equation
$(\epsilon\partial_{tt}^2-\Delta _y )\widetilde \Psi_\epsilon=0$ with initial data $(\Psi_{0,\epsilon},\Psi_{1,\epsilon})$.
Accordingly, we are going to study the integral
\[\begin{array}{l}
\ds\int_0^\infty\ds\int_{\mathbb R^d\times \mathbb R^d} f_\epsilon\nabla_x\Phi_{0,\epsilon} \chi(t,x,v)\ud v\ud x\ud t
\\[.4cm]\qquad=
\ds\int_0^\infty\ds\int_{\mathbb R^d} R_\epsilon(t,x)(\nabla_x\Phi_{0,\epsilon}^{\mathrm{main}} +\nabla_x\Phi_{0,\epsilon}^{\mathrm{rem}})(t,x)\ud x\ud t
\end{array}\]
with $$R_\epsilon(t,x)=\ds\int_{\mathbb R^d} f_\epsilon \chi(t,x,v)\ud v$$
where $\chi$ is a given trial function, supported in $\{0\leq t\leq M,\ |x|\leq M,\ |v|\leq M\}$ for some $0<M<\infty$. 

We observe that 
 \[
  \nabla_x\Big(C_3\ds\frac{\theta_\epsilon-1}{|\cdot|^2}\ast g\Big)=\Big(
  \ds\frac{\nabla_x\theta_\epsilon}{|\cdot|^2}-2(\theta_\epsilon-1)\ds\frac{\cdot}{|\cdot|^4}\Big)\ast g.
 \]
 Thus, by using  \eqref{sigma1} with $d=3$ and $p=2$, we are led to 
 \[
|\nabla_x\Phi_{0,\epsilon}^{\mathrm{rem}}(t,x)|
\leq C \epsilon^{3/4} \left(\ds\int_{\mathbb R^d}\Big| \left(\delta_\epsilon\ast \ds\int_{\mathbb R^n} \sigma_2(y)\widetilde \Psi_\epsilon(t,\cdot,y)\ud y\right)(x')\Big|^2\ud x'\right)^{1/2}
.
\]
However, by \eqref{Psi0} we have
\[\begin{array}{l}
\left\|\delta_\epsilon\ast \ds\int_{\mathbb R^n} \widetilde \Psi_\epsilon\sigma_2(y)\ud y \right\|_{L^\infty([0,\infty);L^2(\R^3))}
\\[.4cm]
\qquad\leq
\|\delta_\epsilon\|_{L^1(\mathbb R^3)}
\|\sigma_2 \|_{L^{2n/(n+2)}(\mathbb R^n)}
\ds\sup_{t\geq 0}\left(\ds\int_{\mathbb R^d}
\| \widetilde \Psi_\epsilon(t,x,\cdot)\|^2_{L^{2n/(n-2)}(\mathbb R^n)}
\ud x\right)^{1/2}
\\[.4cm]
\qquad\leq 
 C \| \sigma_2 \|_{L^{(n+2)/2n}(\R^n)}   \big( \bar{ \mathcal E}_{0}^{\mathrm{vib}} \big)^{1/2}.\end{array}\]
 It implies that $\nabla_x\Phi_{0,\epsilon}^{\mathrm{rem}}(t,x)$ converges uniformly on $(0,\infty)\times \mathbb R^d$ to 0.
 Since $R_\epsilon$ is clearly bounded in $L^1((0,\infty)\times \mathbb R^d\times\mathbb R^d)$, 
 we conclude that $$\ds\int_0^\infty\ds\int_{\mathbb R^d} R_\epsilon \nabla_x\Phi_{0,\epsilon}^{\mathrm{rem}}\ud x\ud t\xrightarrow[\epsilon\to 0]{} 0.$$
 
 We need a more refined estimate to deal with the leading term $\Phi_{0,\epsilon}^{\mathrm{main}}$.
 We begin with
 \[\begin{array}{l}
 \Big|\ds\int_0^\infty\ds\int_{\mathbb R^d} R_\epsilon \nabla_x\Phi_{0,\epsilon}^{\mathrm{main}}\ud x\ud t\Big|
\\[.4cm]
\qquad \leq
 \left(\ds\int_{\mathbb R^d} \left(\ds\int_0^\infty |R_\epsilon|^{p'}\ud t\right)^{2/p'}\ud x\right)^{1/2}
 \left(\ds\int_{\mathbb R^d} \left(\ds\int_0^\infty |\nabla_x\Phi_{0,\epsilon}^{\mathrm{main}}|^p\ud t\right)^{2/p}\ud x\right)^{1/2}.
 \end{array}\]
 We realize that the components of  $\nabla_x\Phi_{0,\epsilon}^{\mathrm{main}}$  are given by the solutions $\Upsilon_{j,\epsilon}$
 of the wave equation
 $$(\epsilon \partial^2_{t}-\Delta_y 
 )\Upsilon_{j,\epsilon}=0$$
 with data
 $$
 \Upsilon_{j,\epsilon}(0,x,y)=\partial_{x_j}\ds\frac{C_3}{|\cdot|^2}\ast \delta_\epsilon\ast \Psi_{0,\epsilon}(x,y),\qquad
 \partial_t \Upsilon_{j,\epsilon}(0,x,y)=\partial_{x_j}\ds\frac{C_3}{|\cdot|^2}\ast \delta_\epsilon\ast \Psi_{1,\epsilon}(x,y),
 $$
 and the space variable $x\in\mathbb R^3$ has only the role of a parameter.
 It satisfies the following Strichartz estimate
 \[
\ds\frac{1}{\epsilon^{1/(2p)}}\left(\ds\int_0^\infty\left(\ds\int_{\mathbb R^n}|\Upsilon_\epsilon(t,x,y)|^q\ud y\right)^{p/q}\ud t\right)^{1/p}\leq  C \sqrt{\mathscr E^{\mathrm{vib}}_{1,\epsilon}(x)}\]
where 
 \[
\mathscr E^{\mathrm{vib}}_{1,\epsilon}(x)=\epsilon\ds\int_{\mathbb R^n}|\partial_t\Upsilon_{\epsilon}(0,x,y)|^2\ud y+ \ds\int_{\mathbb R^n}|\nabla_y\Upsilon_{\epsilon}(0,x,y)|^2\ud y
\]
(for admissible exponents as detailed above).
The Fourier transform of 
 $x\mapsto \nabla_x
  \frac{C_3}{|x|^2}$ is  $\frac{\xi}{|\xi|} $, 
    see \cite[Th. 5.9]{LL}, which implies that the convolution operator
$g\mapsto  \nabla_x
  \frac{C_3}{|x|^2}\ast g$, is an isometry  from $L^2(\R^3)$ to $(L^2(\R^3))^3$.
  Furthermore, we  have $\|\delta_\epsilon\ast g\|_{L^2(\mathbb R^3)}\leq \|\delta_\epsilon\|_{L^1(\mathbb R^3)}\| g\|_{L^2(\mathbb R^3)}=\| g\|_{L^2(\mathbb R^3)}$.
 It follows  that
$$\| \nabla_y  \Upsilon_{\epsilon}(0) \|_{L^2(\R^3_x\times \R^n_y)} \leq \| \nabla_y \Psi_{0,\epsilon} \|_{L^2(\R^3_x\times \R^n_y)},\qquad
\| \partial_t  \Upsilon_{\epsilon}(0) \|_{L^2(\R^3_x\times \R^n_y)} \leq \|  \Psi_{1,\epsilon} \|_{L^2(\R^3_x\times \R^n_y)}.$$ 
 Strichartz' estimate
 then leads to
 $$ \left(\ds\int_{\mathbb R^d} \left(\ds\int_0^\infty |\nabla_x\Phi_{0,\epsilon}^{\mathrm{main}}|^p\ud t\right)^{2/p}\ud x\right)^{1/2}
 \leq C \epsilon^{1/(2p)}\sqrt{\mathscr E^{\mathrm{vib}}_{0,\epsilon}}\leq C \epsilon^{1/(2p)}\sqrt{\bar{\mathscr E}^{\mathrm{vib}}_{0}}.
 $$
 Since  $f_\epsilon$ is bounded in $L^\infty(0,\infty;L^p(\mathbb R^d\times\mathbb R^d))$ for all $1 \leq p \leq \infty$,
 and $\chi$ is bounded and compactly supported we conclude that
 $$\ds\int_0^\infty\ds\int_{\mathbb R^d} R_\epsilon \nabla_x\Phi_{0,\epsilon}^{\mathrm{main}}\ud x\ud t\xrightarrow[\epsilon\to0]{}0.$$
 (Note that the same argument can be applied to show that $\nabla_x\Phi_{0,\epsilon}^{\mathrm{rem}}$ vanishes faster than what has been obtained with the mere energy estimate.) 
\end{Proof}

Next, we study the non linear acceleration term.
Let us set 
$$\widetilde\rho_\epsilon(t,x)=\delta_\epsilon\ast\delta_\epsilon\ast\ds\int_0^{t/\sqrt\epsilon}\rho_\epsilon(t-s\sqrt\epsilon,x)\ q(s)\ud s.$$
It is clear, with Lemma \ref{qL1},  that $\widetilde\rho_\epsilon$ inherits from $\rho_\epsilon$ the uniform estimate  $L^\infty(0,\infty;L^1\cap L^{5/3}(\mathbb R^3))$.
We also denote $E(x)=\frac1{4\pi}\frac{1}{|x|}$, the elementary solution of the  operator $-\Delta_x$ in $\mathbb R^3$.
Note that  $\nabla_x E(x)=-\frac{x}{4\pi|x|^3} $.
Bearing in mind  Lemma~\ref{kernel}, the self--consistent field can be split as follows
\begin{equation}\label{decomp}
\ds\frac{1}{\epsilon}\nabla_x\mathcal L_\epsilon(f_\epsilon)(t,x)=
 \Big[ \nabla_x \Big(\ds\frac{C_3\theta_\epsilon}{|\cdot|^2}\ast\ds\frac{C_3\theta_\epsilon}{|\cdot|^2} \Big)- \nabla_x E\Big]\ast\widetilde\rho_\epsilon(t,x)
 + \nabla_x E\ast\widetilde\rho_\epsilon(t,x).
\end{equation}
In the right hand side, the $L^r$ norm of the first term is dominated by $\|\widetilde\rho_\epsilon\|_{L^\infty([0,\infty;L^1(\mathbb R^3))}\big\|\big[...\big]\big \|_{L^r(\mathbb R^3)}$, hence, owing to Lemma 
Lemma~\ref{kernel} it tends to 0 as $\epsilon\to 0$ in $L^\infty(0,\infty;L^r(\mathbb R^3))$ for any $3/2<r<\infty$. 
Next, Lemma \ref{HLS} tells us that
\[\textrm{
$\nabla_x E\ast\widetilde\rho_\epsilon $ is bounded in $L^\infty ( 0,\infty;L^{15/4}(\mathbb R^3))$}.
\]
Therefore, adapting the reasoning made in the previous sections, we deduce that we can extract a subsequence, such that, for any trial function $\chi\in L^{p'}(\mathbb R^3\times\mathbb R^3)$,
$1/p'+1/p=1$, $1< p< \infty$,
\[
\ds\lim_{\epsilon\to 0}\ds\int_{\mathbb R^3\times\mathbb R^3}f_\epsilon(t,x,v)\chi(x,v)\ud v\ud x=
\ds\int_{\mathbb R^3\times\mathbb R^3}f(t,x,v)\chi(x,v)\ud v\ud x\]
holds uniformly on $[0,T]$, for any $0\leq T<\infty$.
Since the uniform estimate on the kinetic energy imply the tightness of $f_\epsilon$ with respect to the velocity variable,  we also have 
\[
\ds\lim_{\epsilon\to 0}\ds\int_{\mathbb R^3}\rho_\epsilon(t,x,v)\zeta(x)\ud x=
\ds\int_{\mathbb R^3}\rho(t,x)\zeta(x)\ud x,\qquad \rho(t,x)=\ds\int_{\mathbb R^3} f(t,x,v)\ud v,\]
 uniformly on $[0,T]$, for any $0\leq T<\infty$ and any $\zeta \in L^q(\mathbb R^3)$, $q\geq 5/2$ or $\zeta\in C_0(\mathbb R^3)$.
 Clearly, for any  $\zeta\in C^\infty_c(\mathbb R^3)$,  $\delta_{\epsilon}  \ast \delta_{\epsilon}  \ast \zeta$ converges to $\zeta$ in $L^q(\mathbb R^3)$, $5/2\leq q <\infty$,
 and in 
 $C_0(\mathbb R^3)$. 
Therefore 
\[\ds\int_{\mathbb R^3} (\delta_{\epsilon}  \ast \delta_{\epsilon}  \ast \rho_\epsilon) (t,x)\zeta(x)\ud x=
\ds\int_{\mathbb R^3}   \rho_\epsilon (t,x)\ (\delta_{\epsilon}  \ast \delta_{\epsilon}  \ast\zeta)(x)\ud x
\xrightarrow[\epsilon\to 0]{} \kappa\ \ds\int_{\mathbb R^3}   \rho (t,x)\ \zeta(x)\ud x
\]
uniformly in $[0,T]$.
Then, we look at the difference
\[\begin{array}{l}
\left|\ds\int_{\mathbb R^3} \widetilde  \rho_\epsilon (t,x)\zeta(x)\ud x- \kappa\ds\int_{\mathbb R^3} \rho (t,x)\zeta(x)\ud x\right|
\\[.4cm]
\qquad
\leq\ds\int_0^{t/\sqrt\epsilon}\left|\ds\int_{\mathbb R^3} (\delta_{\epsilon}  \ast \delta_{\epsilon}  \ast  \rho_\epsilon) (t-\sqrt\epsilon s,x)\zeta(x)\ud x
-\ds\int_{\mathbb R^3}   \rho (t-\sqrt\epsilon s,x)\zeta(x)\ud x
\right|\  |q(s)|\ud s
\\[.4cm]
\qquad\qquad
+\ds\int_0^{t/\sqrt\epsilon}\left| \ds\int_{\mathbb R^3} \rho (t-\sqrt\epsilon s,x)\zeta(x)\ud x-  \ds\int_{\mathbb R^3} \rho (t,x)\zeta(x)\ud x\right|\ |q(s)|\ud s
\\[.4cm]
\qquad\qquad
+\ds\int_{t/\sqrt\epsilon}^\infty |q(s)| \ud s\ \left|\ds\int_{\mathbb R^3} \rho (t,x)\zeta(x)\ud x\right|.
\end{array}\]
Let us denote by $\mathrm {I}_\epsilon(t)$, $\mathrm {II}_\epsilon(t) $ and $\mathrm {III}_\epsilon(t)$ the three integrals in the right hand side.
By using Lemma \ref{qL1} and the available estimates, we 
obtain, for any $0\leq t\leq T<\infty$ 
\[|\mathrm {I}_\epsilon(t)| \leq \|q\|_{L^1([0,\infty))}\ds\sup_{0\leq u\leq T}\left|\ds\int_{\mathbb R^3}\big( \delta_{\epsilon}  \ast \delta_{\epsilon}  \ast  \rho_\epsilon-\rho\big)(u,x)\zeta(x)\ud x\right|
\xrightarrow[\epsilon\to 0]{}0,\]
while a direct application of the Lebesgue theorem shows that, for any  $0< t\leq T<\infty$
\[
\ds\lim_{\epsilon\to 0}\mathrm {II}_\epsilon(t) =0=\ds\lim_{\epsilon\to 0} \mathrm {III}_\epsilon(t).
\]
Therefore, for any $\zeta\in L^q(\mathbb R^3)$, $5/2\leq q <\infty$
 and any $\zeta\in C_0(\mathbb R^3)$,
 \[
\ds\lim_{\epsilon\to 0} \ds\int_{\mathbb R^3} \widetilde  \rho_\epsilon (t,x)\zeta(x)\ud x= \kappa\ds\int_{\mathbb R^3} \rho (t,x)\zeta(x)\ud x\]
 holds for a. e. $t\in (0,T)$, with the domination 
 \[\left|\ds\int_{\mathbb R^3} \widetilde  \rho_\epsilon (t,x)\zeta(x)\ud x\right|
 \leq \|\zeta\|_{L^{p'}(\mathbb R^3)}\ds\sup_{\epsilon>0,\ 0\leq t\leq T}\|\rho_\epsilon(t,\cdot)\|_{L^p(\mathbb R^3)} ,\]
 for any $1\leq p\leq 5/3$.

In oder to justify that the limit $f$ is a solution of the Vlasov--Poisson equation, the only difficulty relies on the treatment of the non linear acceleration term:
\[ \mathrm{NL}_\epsilon(\chi) =\int_0^\infty \int_{\mathbb R^3 \times \mathbb R^3} f_{\epsilon} \nabla_x \ds\frac1\epsilon\mathcal{L}_\epsilon(f_{\epsilon})\cdot\nabla_v \chi \ud v \ud x\ud t\]
where $\chi$ is a trial function in $\chi\in C^\infty_c([0,\infty)\times \mathbb R^d \times \mathbb R^d)$.
Bearing in mind \eqref{decomp}, it is convenient to rewrite
\[\begin{array}{l} \mathrm{NL}_\epsilon(\chi) = \ds\int_0^\infty\ds\int_{\mathbb R^3 } 
\left(\int_{\mathbb R^3 }
f_{\epsilon}\nabla_v \chi \ud v\right)\cdot \nabla_x E\ast \widetilde \rho_\epsilon  \ud x\ud t+ \mathscr R_\epsilon,
\qquad
\ds\lim_{\epsilon\to 0}  \mathscr R_\epsilon=0.
\end{array}\]
Lemma \ref{HLS} implies that $\nabla_x E\ast \widetilde \rho_\epsilon$ is bounded in $L^\infty(0,T;L^{15/4}(\mathbb R^3))$.
For  $\mu>0$, we introduce the  cut--off function $\widetilde \theta_\mu(x)=\theta(x/\mu)$. 
Then we split 
\[
\nabla_x E\ast \widetilde \rho_\epsilon (t,x) =\ds\int_{\mathbb R^3} \widetilde \theta_\mu(x-y)\ds\frac{x-y}{4\pi|x-y|^3}\widetilde \rho_\epsilon(t,y)\ud y
+
\ds\int_{\mathbb R^3}\big(1- \widetilde \theta_\mu(x-y)\big)\ds\frac{x-y}{4\pi|x-y|^3}\widetilde \rho_\epsilon(t,y)\ud y
.\] 
The first term in the right hand side can be made arbitrarily small in $L^p$ norm, $1\leq p\leq 5/3$, uniformly with respect to $\epsilon$, since
it can be dominated by 
\[
\left\| \ds\int_{|x-y|\leq 2\mu} \ds\frac{x-y}{4\pi|x-y|^3}\widetilde \rho_\epsilon(t,y)\ud y\right\|_{L^p(\mathbb R^3)}
\leq \|\widetilde \rho_\epsilon(t,\cdot)\|_{L^p(\mathbb R^3)}\ \ds\int_{|x-y|\leq 2\mu }\ds\frac{\ud y}{4\pi |x-y|^2}
\leq C\ \mu.
\]
In the second term, for fixed $x\in\mathbb R^3$ and $\mu$, $y\mapsto \big(1-\widetilde \theta_\mu(x-y)\big)\frac{x-y}{4\pi|x-y|^3|}\mathbf 1_{|x-y|\geq \mu}$ is a continuous  function
which vanishes as $|y|\to \infty$, 
so that,  for any $t>0$,
$$\lim_{\epsilon\to 0} \ds\int_{\mathbb R^3}\big(1-\widetilde \theta_\mu(x-y)\big) \ds\frac{x-y}{4\pi|x-y|^3}\widetilde \rho_\epsilon(t,y)\ud y
=\ds\int_{\mathbb R^3}\big(1-\widetilde \theta_\mu(x-y)\big)  \ds\frac{x-y}{4\pi|x-y|^3}\rho(t,y)\ud y.$$
By standard arguments of integration theory (see for instance \cite[Th. 7.61]{Th}), we deduce that (a suitable subsequence of) 
$\nabla_x E\ast \widetilde \rho_\epsilon$ converges to $\nabla_xE\ast\rho$ a. e. and strongly 
in $L^p_{\mathrm{loc}}((0,T)\times\mathbb R^3)$, for any $1\leq p<15/4$.
On the other hand, $\int_{\mathbb R^3 }
f_{\epsilon}\nabla_v \chi \ud v$ is compactly supported and converges  
to $\int_{\mathbb R^3 }
f_{\epsilon}\nabla_v \chi \ud v$ weakly in any $L^q((0,T)\times \mathbb R^3)$.
(In fact this convergence, as well as $\rho_\epsilon\rightarrow \rho$ can be shown to hold strongly, by applying average lemma techniques, see \cite[Th. 5]{DPLM}.)
We conclude that 
\[\ds\lim_{\epsilon\to 0}\mathrm{NL}_\epsilon(\chi) =
 \int_0^\infty\int_{\mathbb R^3 } 
\left(\int_{\mathbb R^3 }
f\nabla_v \chi \ud v\right)\cdot \nabla_x E\ast  \rho  \ud x\ud t.\]
It ends the proof of Theorem \ref{cvth2}.\end{ProofOf}

\bibliography{VW}
\bibliographystyle{plain}

\end{document}